# A new six parameters generalized Weibull-Kumaraswamy-{Lomax) Distribution with a Novel Linear Series Expansion of Its Probability Density

Iman M. Attia *

imanattiathesis1972@gmail.com

**Department of Mathematical Statistics, Faculty of Graduate Studies for Statistical Research, Cairo University, Egypt*

## Abstract

One method of generalizing a distribution is to add new parameter to the baseline distribution. In this paper the author discusses this methodology to create a new six parameters generalized Weibull-Kumaraswamy-{Lomax) Distribution. The author utilizes the approach of the T-X{Y} family to add new parameters the baseline distribution Kumaraswamy distribution defined on the unit interval with two shape parameters $a$ & $b > 0$. The author used the Lomax distribution with two parameters $\alpha$ & $\lambda > 0$ , scale and shape parameters respectively, to link the baseline distribution with the generator distribution. The generator distribution is the Weibull distribution with the two parameters, $c$ & $\gamma > 0$ , shape and scale parameters respectively as will be explained in the paper. The result is a new distribution with six parameters. The author derived the PDF, CDF, and quantile functions. The author used a new technique to expand the pdf in a linear series expansion that is simple and does not require special function. This makes the computational statistics for the raw moments easier and simpler. The author discusses the Shannon entropy and the paper contains figures for the PDF, CDF and the quantile functions.

## Introduction

Many authors developed generalization by introducing families like the beta generated family of distribution BG-distributions which is considered as a generalization of the distribution of the order statistics to mention some of these distributions are Beta-normal (Eugene et al., 2002), Beta-Gumbel (Nadarajah & Kotz, 2004), Beta Frechet (Barreto-Souza et al., 2008), Beta-Exponential (Nadarajah & Kotz, 2006), Beta-Weibull (Lee et al., 2007), Beta-Pareto (Akinsete et al., 2008), Beta-generalized exponential (Barreto-Souza et al., 2010), Beta-modified Weibull (Silva et al., 2010), and Beta-cauchy (Alshawarbeh et al., 2013). Also the Kumaraswamy generated distributions, Kw-G distributions are distributions that tilize the Kumaraswamy distribution to generate new distributions. Some of them are Kumaraswamy-normal (Cordeiro & De Castro, 2011), Kumaraswamy-gamma (Cordeiro & De Castro, 2011), Kukaraswamy-Pareto (Bourguignon et al., 2013), and Kumaraswamy-Gumbel (Cordeiro et al., 2012). While the BG- distributions and KwG-generated distributions use the Beta and Kumaraswamy distributions as the generator and both are defined on the unit interval, (Alzaatreh et al., 2013) introduced general method allowing to use any continuous PDF as a generator. Let $r(t)$ be a PDF of a random variable $T$ defined on the interval $[a,b]$ which is $-\infty < a < b < +\infty$. Let $W(F(x))$ be a function of $F(x)$ of any random variable $X$ so that it satisfies the conditions $W(F(x)) \in [a,b]$ and $W(F(x))$ is differential and monotonically non decreasing and $W(F(x)) \to a$ as $x \to -\infty$ and $W(F(x)) \to b$ as $x \to +\infty$ so the CDF of a new family of distribution can be defines as $G(x) = \int_a^{W(F(x))} r(t)\,dt = R\left(W\big(F(x)\big)\right)$, and the corresponding PDF is $g(x) = \left\{\frac{d}{dx} W\big(F(x)\big)\right\} r\{W\big(F(x)\big)\}$. This family was named as transformed-transformer, where the $W\big(F(x)\big)$ is the transformer that transformed the $r(t)$ into a new CDF, $G(x)$. This $W\big(F(x)\big)$ was suggested to be standard quantile function for some well-known standard distributions like the standard exponential, standard Weibull, standard log-logistic, and the standard Dagum. Also (Lee et al., 2013) used exponentiated CDF of x, $F^{\alpha}(x)$ , in these quantiles. (Aljarrah et al., 2014) extended this approach by applying $W = Q_Y(F(x))$ , where Q is a quantile of the Y variable that act like a link between the X baseline distribution and the T generator distribution. Here the addition of the new parameters come from considering the Y distribution as any continuous distribution with specific parameters and not a standard distribution. This approach yields the following CDF and PDF

$$G_{X_g}(x) = \int_a^{Q_Y(F(x))} r(t)\ dt \ = \ R\left(Q_Y\big(F_X(x)\big)\right) \ , \qquad x \in support\ of\ the\ X\ distribution$$

$g_{X_g}(x) = f_X(x)\, Q_Y'(p)\, f_T(t)\,;\ t = Q_Y\big(F(x)\big)\ \&\ p = F_X(x)$ , the PDF can also be written as

$$g_{X_g}(x) = f_X(x) \frac{r\left(Q_Y\big(F_X(x)\big)\right)}{f_Y\left(Q_Y\big(F_X(x)\big)\right)}\ ;\ \ r\left(Q_Y\big(F_X(x)\big)\right) = \ f_T\left(Q_Y\big(F_X(x)\big)\right)$$

$$Q_Y\big(F_Y(x)\big) = x\,,\ \because\ Q_Y'\big(F_Y(x)\big) f_Y(x) = \frac{d}{dx}\ (x)\ ,\therefore\ Q_Y'\big(F_Y(x)\big) f_Y(x) = 1, \&\ Q_Y'\big(F_Y(x)\big) = \frac{1}{f_Y(x)}$$

(Osatohanmwen et al., 2021) discussed the generalized family of T-kumaraswamy family of distributions with the following distributions: Weibull-Kumaraswamy{exponential}, log-logistic-Kumaraswamy{exponential},exponential-kumaraswamy{log-logistic}, normal-kumaraswamy{logistic}, and logistic-Kumaraswamy{extreme-value}. (Sudsila et al., 2022) discussed the generalized family of distribution based on T-ToppLeone family of distributions and the discussed distributions were: Weibull-TL-{exponential}, log-logistic-TL{exponential}, Logistic-TL{extreme-value}, exponential-TL{log-logistic}, and normal-TL{logistic}.

Many authors tackled the generalization of distributions on unit interval like:(McDonald, 2008) discussed the generalized beta distribution of the first kind, (Carrasco et al., 2010) discussed generalized Kumaraswamy distribution, (El-Sherpieny & Ahmed, 2014) discussed Kumaraswamy-Kumaraswamy distribution, (Elgarhy et al., 2018) discussed exponentiated generalized Kumaraswamy distribution. Some of the major importance of the unit distributions are utilizing these distributions as conjugate prior to some well-known discrete distribution in Bayesian statistics and inference, and fitting different shapes and models of the percentages, proportions, fractions and ratios with random behavior.

To the best of the author knowledge the Weibull-Kumaraswamy{Lomax} was not discussed before in the literature. Also, the new approach of the linear series expansion introduced by the author was not tackled before.

This paper is structured in 7 sections. Section 1 elucidates the methodology of derivation, high-lightening the PDF, CDF, quantile and hazard functions. section 2 illustrates the linear series expansion and utilizing this expansion to obtain the raw moments. Section 3 explains the Shannon entropy. Section 4 depicts some figures for the PDF, CDF and hazard functions for this new WKL distribution. Section 5 enlightens the results of real data analysis and discussing these data in context of the new distribution and in comparison with the baseline distribution, Beta distribution and GOMBUR type I distribution. Section 6 comprehends the conclusion. Section 7 states future works.

# 1. Methodology (Methods and Materials)

## Derivation of Weibull-Kumaraswamy {Lomax} Distribution

The generator distribution is the Weibull (the distribution of the T variable), the baseline distribution is the Kumaraswamy (the distribution of the X variable) and the distribution which can be considered as a link between the T and X distributions is the Lomax distribution (the distribution of the Y variable, the distribution of the quantile function).

The **Weibull Distribution** has the following PDF, CDF and quantile functions $f_T(t) = \frac{c}{\gamma}\left(\frac{t}{\gamma}\right)^{c-1} e^{-\left(\frac{t}{\gamma}\right)^c}$ , $F_T(t) = 1 - e^{-\left(\frac{t}{\gamma}\right)^c}$ , $Q_T(p) = \gamma(-\ln(1-p))^{\frac{1}{c}}$ , respectively, where; $t > 0$, $0 < p < 1$, the shape and scale parameters are $c \;\&\; \gamma > 0$ respectively.

The **Kumaraswamy Distribution** has the following PDF, CDF, and quantile functions respectively: $f_X(x) = abx^{a-1}(1-x^a)^{b-1}$ , $F_X(x) = 1-(1-x^a)^b$ , $Q_X(p) = \left(1-[1-p]^{\frac{1}{b}}\right)^{\frac{1}{a}}$ ,respectively, where $0 < x < 1$, $0 < p < 1$, the two shape parameters are $a \;\&\; b > 0$.

The **Lomax Distribution** has the following PDF, CDF, and quantile functions $f_Y(y) = \frac{\alpha}{\lambda}\left[1+\frac{y}{\lambda}\right]^{-(\alpha+1)}$ , $F_Y(y) = 1-\left[1+\frac{y}{\lambda}\right]^{-\alpha}$ , $Q_Y(p) = \lambda\left[(1-p)^{\frac{-1}{\alpha}}-1\right]$ , respectively, where; $y > 0$, $0 < p < 1$, the shape and scale parameters are $\alpha \;\&\; \lambda > 0$ respectively

Use the approach discussed by Aljarrah 2014 as shown in the following steps to derive the PDF.

**Step 1**: substitute the CDF of Kumaraswamay into the quantile function of the link distribution (quantile function of Y), to obtain the t variable.

$$T = Q_Y[F_X(x)] = \lambda\left[(1-p)^{\frac{-1}{\alpha}}-1\right] = \lambda\left[\left(1-F_X(x)\right)^{\frac{-1}{\alpha}}-1\right]$$

$$= \lambda\left[\left(1-\{1-(1-x^a)^b\}\right)^{\frac{-1}{\alpha}}-1\right] = \lambda\left[\left((1-x^a)^b\right)^{\frac{-1}{\alpha}}-1\right] = \lambda\left[(1-x^a)^{\frac{-b}{\alpha}}-1\right]$$

**Step 2:** substitute this T variable into the PDF of the Weibull distribution

$$f_T(t) = \frac{c}{\gamma}\left(\frac{t}{\gamma}\right)^{c-1} e^{-\left(\frac{t}{\gamma}\right)^c}$$

$$f_T(x) = \frac{c}{\gamma}\left(\frac{\lambda\left[(1-x^a)^{\frac{-b}{\alpha}}-1\right]}{\gamma}\right)^{c-1} exp\left[-\left(\frac{\lambda\left[(1-x^a)^{\frac{-b}{\alpha}}-1\right]}{\gamma}\right)^c\right]$$

**Step 3**: get the derivative of the quantile function which is the reciprocal of the PDF of the quantile distribution:

$$\frac{d}{dx}\left(Q_Y(p)\right) = \frac{d}{dx}\left(Q_Y\left(F_X(x)\right)\right) = \frac{1}{f_Y(t)} = \frac{1}{\frac{\alpha}{\lambda}\left[1+\frac{\lambda\left[(1-x^a)^{\frac{-b}{\alpha}}-1\right]}{\lambda}\right]^{-(\alpha+1)}}$$

$$= \frac{\lambda}{\alpha\left[1+(1-x^a)^{\frac{-b}{\alpha}}-1\right]^{-(\alpha+1)}} = \frac{\lambda}{\alpha[(1-x^a)]^{\frac{b(\alpha+1)}{\alpha}}} = \frac{\lambda}{\alpha[1-x^a]^{b+\frac{b}{\alpha}}}$$

**Step 4**: substitute into the equation of $g(x)$ :

$$g(x) = f_X(x)\frac{r\left(Q_Y\left(F_X(x)\right)\right)}{f_Y\left(Q_Y\left(F_X(x)\right)\right)} = f_X(x)\frac{r\left(Q_Y\left(F_X(x)\right)\right)}{f_Y(t)}\ ,\ \text{where } r\left(Q_Y\left(F_X(x)\right)\right) = f_T(t)$$

$$g(x) = abx^{a-1}(1-x^a)^{b-1}\frac{\frac{c}{\gamma}\left(\frac{\lambda\left[(1-x^a)^{\frac{-b}{\alpha}}-1\right]}{\gamma}\right)^{c-1} exp\left[-\left(\frac{\lambda\left[(1-x^a)^{\frac{-b}{\alpha}}-1\right]}{\gamma}\right)^c\right]}{\frac{\alpha}{\lambda}[1-x^a]^{b+\frac{b}{\alpha}}}$$

After rearrangement and for convenient the $x_g = x$ in the following PDF and CDF where x is defined on the unit interval:

$\boldsymbol{g_{X_g}(x_g;\ a,b,c,\alpha,\gamma,\lambda)} =$

$$\frac{abc\lambda}{\alpha\gamma}x^{a-1}(1-x^a)^{-1-\frac{b}{\alpha}}\left(\frac{\lambda\left[(1-x^a)^{\frac{-b}{\alpha}}-1\right]}{\gamma}\right)^{c-1} exp\left[-\left(\frac{\lambda\left[(1-x^a)^{\frac{-b}{\alpha}}-1\right]}{\gamma}\right)^c\right]$$

$$\boldsymbol{G_{X_g}(x_g;\ a,b,c,\alpha,\gamma,\lambda)} = F_T\ (t) = 1 - exp\left[-\left(\frac{\lambda\left[(1-x^a)^{\frac{-b}{\alpha}}-1\right]}{\gamma}\right)^c\right]; 0 < x < 1$$

The quantile function is:

$$\boldsymbol{Q}_{\boldsymbol{X}_{\boldsymbol{g}}}(\boldsymbol{u};\boldsymbol{a},\boldsymbol{b},\boldsymbol{c},\boldsymbol{\alpha},\boldsymbol{\gamma},\boldsymbol{\lambda}) = \left[1-\left\{1+\frac{\gamma}{\lambda}[-\ln(1-u)]^{\frac{1}{c}}\right\}^{-\frac{\alpha}{b}}\right]^{\frac{1}{a}};\ 0<u<1, a,b,c,\alpha,\gamma,\lambda>0$$

The hazard function is

$$\boldsymbol{h}(\boldsymbol{X}_{\boldsymbol{g}}) = \frac{\boldsymbol{g}_{\boldsymbol{X}_{\boldsymbol{g}}}}{1-\boldsymbol{G}_{\boldsymbol{X}_{\boldsymbol{g}}}} = \frac{abc\lambda}{\alpha\gamma}x^{a-1}(1-x^a)^{-1-\frac{b}{\alpha}}\left(\frac{\lambda\left[(1-x^a)^{\frac{-b}{\alpha}}-1\right]}{\gamma}\right)^{c-1}$$

## 2.linear series expansion

**Theorem 1**: linear series expansion of this pdf, is applied to the exponent term and each of the power terms.

**Proof**: using the well-known series expansion for the exponent part:

$$e^{-z} = \sum_{k=0}^{\infty}\frac{(-1)^k}{k!}z^k\ ;\quad call\ \ z = \left[\left(\frac{\lambda(1-x^a)^{\frac{-b}{\alpha}}\left[1-(1-x^a)^{\frac{b}{\alpha}}\right]}{\gamma}\right)^c\right];$$

$$exp\left[\left(\frac{\lambda(1-x^a)^{\frac{-b}{\alpha}}\left[1-(1-x^a)^{\frac{b}{\alpha}}\right]}{\gamma}\right)^c\right] = e^{-z}$$

$$e^{-z} = \sum_{k=0}^{\infty}\frac{(-1)^k}{k!}\left[\left(\frac{\lambda(1-x^a)^{\frac{-b}{\alpha}}\left[1-(1-x^a)^{\frac{b}{\alpha}}\right]}{\gamma}\right)^c\right]^k$$

$$= \sum_{k=0}^{\infty}\frac{(-1)^k}{k!}\frac{\lambda^{ck}(1-x^a)^{\frac{-ckb}{\alpha}}\left[1-(1-x^a)^{\frac{b}{\alpha}}\right]^{ck}}{\gamma^{ck}}$$

$$\sum_{k=0}^{\infty}\frac{(-1)^k}{k!}\left(\frac{\lambda}{\gamma}\right)^{ck}\sum_{h=0}^{\infty}(-1)^h\binom{\frac{-kcb}{\alpha}}{h}[x^a]^h\left[1-(1-x^a)^{\frac{b}{\alpha}}\right]^{ck} =$$

$$\sum_{k=0}^{\infty}\frac{(-1)^k}{k!}\left(\frac{\lambda}{\gamma}\right)^{ck}\sum_{h=0}^{\infty}(-1)^h\binom{\frac{-kcb}{\alpha}}{h}[x^a]^h\sum_{n=0}^{\infty}(-1)^n\binom{kc}{n}\left[(1-x^a)^{\frac{b}{\alpha}}\right]^n =$$

$$\sum_{k=0}^{\infty}\frac{(-1)^k}{k!}\left(\frac{\lambda}{\gamma}\right)^{ck}\sum_{h=0}^{\infty}(-1)^h\binom{\frac{-kcb}{\alpha}}{h}[x^a]^h\sum_{n=0}^{\infty}(-1)^n\binom{kc}{n}\sum_{l=0}^{\infty}(-1)^l\binom{n}{l}\left[(x^a)^{\frac{b}{\alpha}}\right]^l$$

$$\sum_{k=0}^{\infty}\sum_{h=0}^{\infty}\sum_{n=0}^{\infty}\sum_{l=0}^{\infty}\frac{(-1)^{k+h+n+l}}{k!}\left(\frac{\lambda}{\gamma}\right)^{ck}\binom{\frac{-kcb}{\alpha}}{h}\binom{ck}{n}\binom{n}{l}x^{ah+\frac{abl}{\alpha}}$$

Second expand the power part using the generalized binomial expansion formula with the Pochhammer falling factorial

$$(1-x^{a})^{-\left(1+\frac{b}{\alpha}\right)} = \sum_{j=0}^{\infty}(-1)^{j}\binom{-\left(1+\frac{b}{\alpha}\right)}{j}x^{aj}$$

Third expand the remaining power part using the generalized binomial expansion formula.

$$\left[\frac{\lambda\left[(1-x^{a})^{\frac{-b}{\alpha}}-1\right]}{\gamma}\right]^{c-1} = \left[\frac{\lambda\left[(1-x^{a})^{\frac{-b}{\alpha}}\right]\left[1-(1-x^{a})^{\frac{b}{\alpha}}\right]}{\gamma}\right]^{c-1}$$

$$= \left(\frac{\lambda}{\gamma}\right)^{c-1}\left[(1-x^{a})^{\frac{-b}{\alpha}}\right]^{c-1}\left[1-(1-x^{a})^{\frac{b}{\alpha}}\right]^{c-1}$$

$$\left(\frac{\lambda}{\gamma}\right)^{c-1}\sum_{m=0}^{\infty}(-1)^{m}\binom{\frac{-b(c-1)}{\alpha}}{m}x^{am}\sum_{f=0}^{\infty}(-1)^{f}\binom{c-1}{f}\left[(1-x^{a})^{\frac{b}{\alpha}}\right]^{f} =$$

$$\left(\frac{\lambda}{\gamma}\right)^{c-1}\sum_{m=0}^{\infty}(-1)^{m}\binom{\frac{-b(c-1)}{\alpha}}{m}x^{am}\sum_{f=0}^{\infty}(-1)^{f}\binom{c-1}{f}\sum_{p=0}^{\infty}(-1)^{p}\binom{\frac{bf}{\alpha}}{p}[x^{a}]^{p} =$$

$$\left(\frac{\lambda}{\gamma}\right)^{c-1}\sum_{m=0}^{\infty}\sum_{f=0}^{\infty}\sum_{p=0}^{\infty}(-1)^{m+f+p}\binom{\frac{-b(c-1)}{\alpha}}{m}\binom{c-1}{f}\binom{\frac{bf}{\alpha}}{p}\left[x^{a(m+p)}\right]^{1}$$

Now put the three parts together $g_{X_g}(x)$

$$= \frac{abc\lambda}{\alpha\gamma}x^{a-1}\sum_{k=0}^{\infty}\sum_{h=0}^{\infty}\sum_{n=0}^{\infty}\sum_{l=0}^{\infty}\frac{(-1)^{k+h+n+l}}{k!}\left(\frac{\lambda}{\gamma}\right)^{ck}\binom{\frac{-kcb}{\alpha}}{h}\binom{ck}{n}\binom{n}{l}x^{ah+\frac{abl}{\alpha}}\sum_{j=0}^{\infty}(-1)^{j}\binom{-\left(1+\frac{b}{\alpha}\right)}{j}x^{aj}\left(\frac{\lambda}{\gamma}\right)^{c-1}$$

$$\sum_{m=0}^{\infty}\sum_{f=0}^{\infty}\sum_{p=0}^{\infty}(-1)^{m+f+p}\binom{\frac{-b(c-1)}{\alpha}}{m}\binom{c-1}{f}\binom{\frac{bf}{\alpha}}{p}x^{a(m+p)}$$

$$= \frac{abc}{\alpha}\left(\frac{\lambda}{\gamma}\right)^{c}x^{a-1+ah++\frac{abl}{\alpha}+aj+am+ap}$$

$$\sum_{k=0}^{\infty}\sum_{h=0}^{\infty}\sum_{n=0}^{\infty}\sum_{l=0}^{\infty}\sum_{j=0}^{\infty}\sum_{m=0}^{\infty}\sum_{f=0}^{\infty}\sum_{p=0}^{\infty}\frac{(-1)^{k+h+n+l+j+m+f+p}}{k!}\left(\frac{\lambda}{\gamma}\right)^{ck}\binom{\frac{-kcb}{\alpha}}{h}\binom{ck}{n}\binom{n}{l}\binom{-\left(1+\frac{b}{\alpha}\right)}{j}\binom{\frac{-b(c-1)}{\alpha}}{m}\binom{c-1}{f}\binom{\frac{bf}{\alpha}}{p}$$

Integrate this series expansion by integrating the x from 0 to 1 because the support of the distribution is the support of the Kumaraswamy so the integral will pass to the x

$$\frac{abc}{\alpha}\left(\frac{\lambda}{\gamma}\right)^{c}$$

$$\sum_{k=0}^{\infty}\sum_{h=0}^{\infty}\sum_{n=0}^{\infty}\sum_{l=0}^{\infty}\sum_{j=0}^{\infty}\sum_{m=0}^{\infty}\sum_{f=0}^{\infty}\sum_{p=0}^{\infty}\frac{(-1)^{k+h+n+l+j+m+f+p}}{k!}\left(\frac{\lambda}{\gamma}\right)^{ck}\binom{\frac{-kcb}{\alpha}}{h}\binom{ck}{n}\binom{n}{l}\binom{-\left(1+\frac{b}{\alpha}\right)}{j}\binom{\frac{-b(c-1)}{\alpha}}{m}\binom{c-1}{f}\binom{\frac{bf}{\alpha}}{p}$$

$$\int_{0}^{1}x^{a-1+ah++\frac{abl}{\alpha}+aj+am+ap}\,dx$$

This gives $\int_0^1 x^{a-1+ah+\frac{abl}{\alpha}+aj+am+ap} = \left.\frac{x^{a+ah+\frac{abl}{\alpha}+aj+am+ap}}{a+ah++\frac{abl}{\alpha}+aj+am+ap}\right|_0^1 = \frac{1}{a\left(1+h+\frac{bl}{\alpha}+j+m+p\right)}$

$$\frac{abc}{\alpha}\left(\frac{\lambda}{\gamma}\right)^c W_{k,h,n,l,j,m,f,p}\int_0^1 x^{a-1+ah++\frac{abl}{\alpha}+aj+am+ap}dx$$

$$=\frac{abc}{\alpha}\left(\frac{\lambda}{\gamma}\right)^c \frac{W_{k,h,n,l,j,m,f,p}}{a\left(1+h+\frac{bl}{\alpha}+j+m+p\right)}$$

$$g(x)=\frac{abc}{\alpha}\left(\frac{\lambda}{\gamma}\right)^c W_{k,h,n,l,j,m,f,p}\,.\,x^{a-1+ah++\frac{abl}{\alpha}+aj+am+ap} \quad 0<x<1\,, \qquad a,b,\alpha,\gamma,\lambda\,,c>0$$

**Theorem 2: the r<sup>th</sup> non-central moments are given by the following formula:**

$$\boldsymbol{E(x^r)=\frac{abc}{\alpha}\left(\frac{\lambda}{\gamma}\right)^c \frac{W_{k,h,n,l,j,m,f,p}}{a\left(1+h+\frac{bl}{\alpha}+j+m+p\right)+r}}$$

**Proof:**

$$\boldsymbol{E(x^r)=\int_0^1 g(x)x^r dx=\frac{abc}{\alpha}\left(\frac{\lambda}{\gamma}\right)^c W_{k,h,n,l,j,m,f,p}\,.\int_0^1 x^r x^{a-1+ah++\frac{abl}{\alpha}+aj+am+ap}dx}$$

$$\boldsymbol{E(x^r)=\frac{abc}{\alpha}\left(\frac{\lambda}{\gamma}\right)^c \frac{W_{k,h,n,l,j,m,f,p}}{a\left(1+h+\frac{bl}{\alpha}+j+m+p\right)+r}}$$

## 3.Shannon entropy

**Theorem 3**: The entropy of a random variable is a measure of the measure of uncertainty (Shannon, 1948). Shannon's entropy of the random variable X with density g is defined as $E\{-\log(g[x])\}$

**Proof:**

$$shannon\ entropy = E\{-\log(g[x])\}$$

$$\log(g[x]) = \ln a + \ln b + \ln c + \ln\lambda - \ln\alpha - \ln\gamma + (a-1)\ln x - \left(\boldsymbol{1}+\frac{\boldsymbol{b}}{\boldsymbol{\alpha}}\right)\ln(\boldsymbol{1}-\boldsymbol{x^a}) + (c-1)\ln\left[\frac{\boldsymbol{\lambda}\left[(\boldsymbol{1}-\boldsymbol{x^a})^{\frac{-\boldsymbol{b}}{\boldsymbol{\alpha}}}-\boldsymbol{1}\right]}{\boldsymbol{\gamma}}\right]-\left[\left(\frac{\boldsymbol{\lambda}\left[(\boldsymbol{1}-\boldsymbol{x^a})^{\frac{-\boldsymbol{b}}{\boldsymbol{\alpha}}}-\boldsymbol{1}\right]}{\boldsymbol{\gamma}}\right)^c\right]$$

$$\log(g[x]) = \ln a + \ln b + \ln c + \ln\lambda - \ln\alpha - \ln\gamma + (a-1)\ln x - \left(1+\frac{b}{\alpha}\right)\ln(1-x^a) + (c-1)\ln\lambda + (c-1)\ln\left[(1-x^a)^{\frac{-b}{\alpha}}-1\right]-(c-1)\ln\gamma - \left[\left(\frac{\lambda\left[(1-x^a)^{\frac{-b}{\alpha}}-1\right]}{\gamma}\right)^c\right]$$

$$\log(g[x]) = \ln a + \ln b + \ln c + c \ln \lambda - \ln \alpha - c \ln \gamma + (a-1)\ln x - \left(1+\frac{b}{\alpha}\right)\ln(1-x^{a}) + (c-1)\ln\left[(1-x^{a})^{\frac{-b}{\alpha}} - 1\right] - \left(\frac{\lambda}{\gamma}\right)^{c}\left((1-x^{a})^{\frac{-b}{\alpha}} - 1\right)^{c}$$

$$-\log(g[x]) = -\ln a - \ln b - \ln c - c \ln \lambda + \ln \alpha + c \ln \gamma - (a-1)\ln x + \left(1+\frac{b}{\alpha}\right)\ln(1-x^{a}) - (c-1)\ln\left[(1-x^{a})^{\frac{-b}{\alpha}} - 1\right] + \left(\frac{\lambda}{\gamma}\right)^{c}\left((1-x^{a})^{\frac{-b}{\alpha}} - 1\right)^{c}$$

Take the expectation of the above equation.

**First part**: $E\left[-(a-1)\ln x\right] = (1-a)\int_0^1 g(x)\ln x \; dx$

**Second part**: $E\left[\left(1+\frac{b}{\alpha}\right)\ln(1-x^{a})\right] = \left(1+\frac{b}{\alpha}\right)\int_0^1 g(x)\ln(1-x^{a}) \; dx$

**Third part**: $E\left[-(c-1)\ln\left[(1-x^{a})^{\frac{-b}{\alpha}} - 1\right]\right] = (1-c)\int_0^1 g(x)\ln\left[(1-x^{a})^{\frac{-b}{\alpha}} - 1\right] dx$

**Fourth part**: $E\left[\left(\frac{\lambda\left[(1-x^{a})^{\frac{-b}{\alpha}} - 1\right]}{\gamma}\right)^{c}\right] = \int_0^1 g(x)\left[\left(\frac{\lambda}{\gamma}\right)^{c}\left((1-x^{a})^{\frac{-b}{\alpha}} - 1\right)^{c}\right] dx$

In each part pass the expectation to the X variable

**First part**: $(1-a)\int_0^1 g(x)\ln x \; dx =$

$$(1-a)\int_0^1 \frac{abc}{\alpha}\left(\frac{\lambda}{\gamma}\right)^{c} W_{k,h,n,l,j,m,f,p}\, x^{a-1+ah+\frac{abl}{\alpha}+aj+am+ap}\ln x \; dx$$

$$= (1-a)\frac{abc}{\alpha}\left(\frac{\lambda}{\gamma}\right)^{c} W_{k,h,n,l,j,m,f,p}\int_0^1 x^{a-1+ah+\frac{abl}{\alpha}+aj+am+ap}\ln x \; dx\,; \qquad \text{integrate by parts}$$

$f(x) = \ln x$ $\qquad$ $f'(x) = \frac{1}{x}$

$h'(x) = x^{a-1+ah+\frac{abl}{\alpha}+aj+am+ap}$ $\qquad$ $h(x) = \frac{x^{a+ah+\frac{abl}{\alpha}+aj+am+ap}}{a+ah+\frac{abl}{\alpha}+aj+am+ap}$

$$I = \left.\frac{x^{a+ah+\frac{abl}{\alpha}+aj+am+ap}}{a+ah+\frac{abl}{\alpha}+aj+am+ap}\,.\,\ln x\right|_0^1 - \int_0^1 \frac{x^{a+ah+\frac{abl}{\alpha}+aj+am+ap}}{a+ah+\frac{abl}{\alpha}+aj+am+ap}\,.\frac{1}{x}\; dx$$

$$I = 0 - \left.\frac{1}{a+ah+\frac{abl}{\alpha}+aj+am+ap}\,.\,\frac{x^{a+ah+\frac{abl}{\alpha}+aj+am+ap}}{a+ah+\frac{abl}{\alpha}+aj+am+ap}\right|_0^1$$

$$= \frac{-1}{a+ah+\frac{abl}{\alpha}+aj+am+ap}\,.\,\frac{1}{a+ah+\frac{abl}{\alpha}+aj+am+ap}$$

$$\therefore\ (1-a)\int_0^1 g(x)\ln x\ dx = \frac{(a-1)\ \frac{abc}{\alpha}\left(\frac{\lambda}{\gamma}\right)^c\ W_{k,h,n,l,j,m,f,p}}{\left(a+ah+\frac{abl}{\alpha}+aj+am+ap\right)^2}$$

**Second part**: $E\left[\left(1+\frac{b}{\alpha}\right)\ \ln(1-x^a)\right] =$

$$\left(1+\frac{b}{\alpha}\right)\int_0^1 \frac{abc}{\alpha}\left(\frac{\lambda}{\gamma}\right)^c\ \ W_{k,h,n,l,j,m,f,p}\ x^{a-1+ah+\frac{abl}{\alpha}+aj+am+ap}\ln(1-x^a)\ dx$$

$$= \left(1+\frac{b}{\alpha}\right)\frac{abc}{\alpha}\left(\frac{\lambda}{\gamma}\right)^c\ W_{k,h,n,l,j,m,f,p}\int_0^1\ x^{a-1+ah+\frac{abl}{\alpha}+aj+am+ap}\ln(1-x^a)\ \ ;$$

integrate by parts

$f(x)=\ln(1-x^a)$ $\qquad\qquad$ $f'(x)=\frac{-ax^{a-1}}{(1-x^a)}$

$h'(x)=x^{a-1+ah+\frac{abl}{\alpha}+aj+am+ap}$ $\qquad\qquad$ $h(x)=\frac{x^{a+ah+\frac{abl}{\alpha}+aj+am+ap}}{a+ah+\frac{abl}{\alpha}+aj+am+ap}$

$$I = \left.\frac{x^{a+ah+\frac{abl}{\alpha}+aj+am+ap}}{a+ah+\frac{abl}{\alpha}+aj+am+ap}\ \ln(1-x^a)\right|_0^1 - \int_0^1\frac{x^{a+ah+\frac{abl}{\alpha}+aj+am+ap}}{a+ah+\frac{abl}{\alpha}+aj+am+ap}\cdot\frac{-ax^{a-1}}{(1-x^a)}\ \ dx$$

$$I = \ 0 - \frac{-a}{a+ah+\frac{abl}{\alpha}+aj+am+ap}\int_0^1 x^{a+ah+\frac{abl}{\alpha}+aj+am+ap}\ x^{a-1}\ \ \ (1-x^a)^{-1}\,dx$$

Expand : $(1-x^a)^{-1} = \sum_{s=0}^{\infty}(x^a)^s = \sum_{s=0}^{\infty}x^{as}$

$$I = \frac{a}{a+ah+\frac{abl}{\alpha}+aj+am+ap}\int_0^1 x^{a+ah+\frac{abl}{\alpha}+aj+am+ap}\ x^{a-1}\ \ \sum_{s=0}^{\infty}x^{as}\ dx$$

$$I = \frac{a}{a+ah+\frac{abl}{\alpha}+aj+am+ap}\int_0^1\ \sum_{s=0}^{\infty}x^{as}\,x^{a+ah+\frac{abl}{\alpha}+aj+am+a}\ \ x^{a-1}\,dx$$

$$I = \frac{a}{a+ah+\frac{abl}{\alpha}+aj+am+ap}\ \sum_{s=0}^{\infty}\int_0^1\ x^{as}\,x^{a+ah+\frac{abl}{\alpha}+aj+am+a}\ \ x^{a-1}\ dx$$

$$\sum_{s=0}^{\infty}\int_0^1\ x^{as}\ x^{a+ah+\frac{abl}{\alpha}+aj+am+ap}\ x^{a-1}\,dx = \left.\sum_{s=0}^{\infty}\frac{x^{as+a+ah+\frac{abl}{\alpha}+aj+am+ap+a}}{as+a+ah+\frac{abl}{\alpha}+aj+am+ap+a}\right|_0^1$$

$$= \ \sum_{s=0}^{\infty}\frac{1}{as+a+ah+\frac{abl}{\alpha}+aj+am+ap+a}$$

$$I = \frac{a}{a\left(1+h+\frac{bl}{\alpha}+j+m+p\right)}\cdot\sum_{s=0}^{\infty}\frac{1}{as+a+ah+\frac{abl}{\alpha}+aj+am+ap+a}$$

$$E\left[\left(1+\frac{b}{\alpha}\right)\ \ln(1-x^a)\right] = \left(1+\frac{b}{\alpha}\right)\int_0^1 W_{k,h,n,l,j,m,f,p}\ x^{a-1+ah+\frac{abl}{\alpha}+aj+am+ap}\ln(1-x^a)\ dx$$

$$= \frac{\left(1+\frac{b}{\alpha}\right)\frac{abc}{\alpha}\left(\frac{\lambda}{\gamma}\right)^{c} W_{k,h,n,l,j,m,f,p}}{\left(1+h+\frac{bl}{\alpha}+j+m+p\right)} \sum_{s=0}^{\infty} \frac{1}{a\left(2+h+\frac{bl}{\alpha}+j+m+p+s\right)}$$

**Third part**: $E\left[-(c-1)\ln\left[(1-x^{a})^{\frac{-b}{\alpha}}-1\right]\right] =$

$$(1-c)\int_0^1 \frac{abc}{\alpha}\left(\frac{\lambda}{\gamma}\right)^{c} W_{k,h,n,l,j,m,f,p}\, x^{a-1+ah+\frac{abl}{\alpha}+aj+am+ap} \ \ln\left[(1-x^{a})^{\frac{-b}{\alpha}}-1\right] dx$$

$$f(x) = \ln\left[(1-x^{a})^{\frac{-b}{\alpha}}-1\right] \qquad f'(x) = \frac{\frac{-b}{\alpha}(1-x^{a})^{\frac{-b}{\alpha}-1}(-ax^{a-1})}{\left[(1-x^{a})^{\frac{-b}{\alpha}}-1\right]} = \frac{\frac{b}{\alpha}(1-x^{a})^{\frac{-b}{\alpha}-1}(ax^{a-1})}{\left[(1-x^{a})^{\frac{-b}{\alpha}}-1\right]}$$

$$h'(x) = x^{a-1+ah+\frac{abl}{\alpha}+aj+am+ap} \qquad h(x) = \frac{x^{a+ah+\frac{abl}{\alpha}+aj+am+ap}}{a+ah+\frac{abl}{\alpha}+aj+am+ap}$$

$$I = \left.\frac{x^{a+ah+\frac{abl}{\alpha}+aj+am+ap}}{a+ah+\frac{abl}{\alpha}+aj+am+ap} \ \ln\left[(1-x^{a})^{\frac{-b}{\alpha}}-1\right]\right|_0^1$$

$$-\int_0^1 \frac{x^{a+ah+\frac{abl}{\alpha}+aj+am+ap}}{a+ah+\frac{abl}{\alpha}+aj+am+ap} \ \frac{\frac{b}{\alpha}(1-x^{a})^{\frac{-b}{\alpha}-1}(ax^{a-1})}{\left[(1-x^{a})^{\frac{-b}{\alpha}}-1\right]} \ dx$$

$$I = \left.\frac{x^{a+ah+\frac{abl}{\alpha}+aj+am+ap}}{a+ah+\frac{abl}{\alpha}+aj+am+ap} \ \ln(1-x^{a})^{\frac{-b}{\alpha}}\left[1-(1-x^{a})^{\frac{b}{\alpha}}\right]\right|_0^1$$

$$-\int_0^1 \frac{x^{a+ah+\frac{abl}{\alpha}+aj+am+ap}}{a+ah+\frac{abl}{\alpha}+aj+am+ap} \ \frac{\frac{b}{\alpha}(1-x^{a})^{\frac{-b}{\alpha}-1}(ax^{a-1})}{\left[(1-x^{a})^{\frac{-b}{\alpha}}-1\right]} \ dx$$

Rewrite $\frac{\frac{b}{\alpha}(1-x^{a})^{\frac{-b}{\alpha}-1}(ax^{a-1})}{\left[(1-x^{a})^{\frac{-b}{\alpha}}-1\right]} = \frac{\frac{b}{\alpha}(1-x^{a})^{\frac{-b}{\alpha}-1}(ax^{a-1})}{(1-x^{a})^{\frac{-b}{\alpha}}\left[1-(1-x^{a})^{\frac{b}{\alpha}}\right]} = \frac{\frac{b}{\alpha}(1-x^{a})^{-1}(ax^{a-1})}{\left[1-(1-x^{a})^{\frac{b}{\alpha}}\right]}$

$$I = 0 - \frac{a\frac{b}{\alpha}}{a\left(1+h+\frac{bl}{\alpha}+j+m+p\right)} \int_0^1 x^{a+ah+\frac{abl}{\alpha}+aj+am+ap} \ (1-x^{a})^{-1}(x^{a-1})$$

$$\left[1-(1-x^{a})^{\frac{b}{\alpha}}\right]^{-1} \ dx$$

Expand:

$(1-x^{a})^{-1} = \sum_{v=0}^{\infty} x^{va}$

$\left[1-(1-x^{a})^{\frac{b}{\alpha}}\right]^{-1} = \sum_{r=0}^{\infty}\left[(1-x^{a})^{\frac{b}{\alpha}}\right]^{r} = \sum_{r=0}^{\infty}(1-x^{a})^{\frac{br}{\alpha}} = \sum_{r=0}^{\infty}\sum_{d=0}^{\infty}(-1)^{d}\binom{\frac{br}{\alpha}}{d}x^{ad}$

$I = \frac{-b}{\alpha\left(1+h+\frac{bl}{\alpha}+j+m+p\right)} \int_0^1 x^{a+ah+\frac{abl}{\alpha}+aj+am+ap} \ (x^{a-1}) \ \sum_{v=0}^{\infty} x^{va} \ \sum_{r=0}^{\infty}\sum_{d=0}^{\infty}(-1)^{d}\binom{\frac{br}{\alpha}}{d}x^{ad} dx$

rearrange the terms and swap the integral and summation

$$I = \frac{-b}{\alpha\left(1+h+\frac{bl}{\alpha}+j+m+p\right)} \int_0^1 x^{a+ah+\frac{abl}{\alpha}+aj+am+ap} (x^{a-1}) x^{va} x^{ad} \sum_{v=0}^{\infty}\sum_{r=0}^{\infty}\sum_{d=0}^{\infty}(-1)^d \binom{\frac{br}{\alpha}}{d}\, dx$$

$$I = \frac{-b}{\alpha\left(1+h+\frac{bl}{\alpha}+j+m+p\right)} \sum_{v=0}^{\infty}\sum_{r=0}^{\infty}\sum_{d=0}^{\infty}(-1)^d \binom{\frac{br}{\alpha}}{d} \int_0^1 x^{a+ah+\frac{abl}{\alpha}+aj+am+ap} (x^{a-1}) x^{va} x^{ad}\, dx$$

$$\int_0^1 x^{a+ah+\frac{abl}{\alpha}+aj+am+ap} (x^{a-1}) x^{va} x^{ad}\, dx = \left.\frac{x^{a+ah+\frac{abl}{\alpha}+aj+am+ap+av+ad}}{a+ah+\frac{abl}{\alpha}+aj+am+ap+av+ad+a}\right|_0^1$$

$$= \frac{1}{a\left(2+h+\frac{bl}{\alpha}+j+m+p+v+d\right)}$$

$$E\left[-(c-1)\ln\left[(1-x^a)^{\frac{-b}{\alpha}}-1\right]\right]$$

$$= \frac{b(c-1)\frac{abc}{\alpha}\left(\frac{\lambda}{\gamma}\right)^c W_{k,h,n,l,j,m,f,p}}{\alpha\left(1+h+\frac{bl}{\alpha}+j+m+p\right)} \sum_{v=0}^{\infty}\sum_{r=0}^{\infty}\sum_{d=0}^{\infty}(-1)^d \binom{\frac{br}{\alpha}}{d} \frac{1}{a\left(2+h+\frac{bl}{\alpha}+j+m+p+v+d\right)}$$

**Fourth part**: $E\left[\left(\frac{\lambda\left[(1-x^a)^{\frac{-b}{\alpha}}-1\right]}{\gamma}\right)^c\right] = \int_0^1 g(x)\left[\left(\frac{\lambda}{\gamma}\right)^c\left((1-x^a)^{\frac{-b}{\alpha}}-1\right)^c\right]\, dx$

$$= \int_0^1 \frac{abc}{\alpha}\left(\frac{\lambda}{\gamma}\right)^c W_{k,h,n,l,j,m,f,p}\, x^{a-1+ah+\frac{abl}{\alpha}+aj+am+ap} \left[\left(\frac{\lambda}{\gamma}\right)^c\left((1-x^a)^{\frac{-b}{\alpha}}-1\right)^c\right] dx$$

$$= \frac{abc}{\alpha}\left(\frac{\lambda}{\gamma}\right)^c W_{k,h,n,l,j,m,f,p} \left(\frac{\lambda}{\gamma}\right)^c \int_0^1 x^{a-1+ah+\frac{abl}{\alpha}+aj+am+ap} \left((1-x^a)^{\frac{-b}{\alpha}}-1\right)^c dx$$

Rewrite: $\left((1-x^a)^{\frac{-b}{\alpha}}-1\right)^c = (1-x^a)^{\frac{-cb}{\alpha}} \left[1-(1-x^a)^{\frac{b}{\alpha}}\right]^c$

$$(1-x^a)^{\frac{-cb}{\alpha}} \left[1-(1-x^a)^{\frac{b}{\alpha}}\right]^c = \sum_{z=0}^{\infty}(-1)^z \binom{\frac{-cb}{\alpha}}{z} x^{az} \sum_{w=0}^{\infty}(-1)^w \binom{c}{w} (1-x^a)^{\frac{wb}{\alpha}}$$

$$= \sum_{z=0}^{\infty}(-1)^z \binom{\frac{-cb}{\alpha}}{z} x^{az} \sum_{w=0}^{\infty}(-1)^w \binom{c}{w} \sum_{u=0}^{\infty}(-1)^u \binom{\frac{wb}{\alpha}}{u} x^{au}$$

$$= \sum_{z=0}^{\infty}\sum_{w=0}^{\infty}\sum_{u=0}^{\infty}(-1)^{z+w+u} \binom{\frac{-cb}{\alpha}}{z} \binom{c}{w} \binom{\frac{wb}{\alpha}}{u}\, x^{az+a}$$

$$E\left[\left(\frac{\lambda\left[(1-x^a)^{\frac{-b}{\alpha}}-1\right]}{\gamma}\right)^c\right]$$

$$= W_{k,h,n,l,j,m,f,p} \left(\frac{\lambda}{\gamma}\right)^{2c} \int_0^1 x^{a-1+ah+\frac{abl}{\alpha}+aj+am+ap} \sum_{z=0}^{\infty}\sum_{w=0}^{\infty}\sum_{u=0}^{\infty}(-1)^{z+w+u}\binom{\frac{-cb}{\alpha}}{z}\binom{c}{w}\binom{\frac{wb}{\alpha}}{u} x^{az+au}dx$$

$$= W_{k,h,n,l,j,m,f,p} \left(\frac{\lambda}{\gamma}\right)^{2c} \int_0^1 \sum_{z=0}^{\infty}\sum_{w=0}^{\infty}\sum_{u=0}^{\infty}(-1)^{z+w+u}\binom{\frac{-cb}{\alpha}}{z}\binom{c}{w}\binom{\frac{wb}{\alpha}}{u} x^{a-1+ah+\frac{abl}{\alpha}+aj+am+ap} x^{az+au}dx$$

Swap the integration and summation

$$E\left[\left(\frac{\lambda\left[(1-x^{a})^{\frac{-b}{\alpha}}-1\right]}{\gamma}\right)^{c}\right]$$

$$= W_{k,h,n,l,j,m,f,p} \left(\frac{\lambda}{\gamma}\right)^{2c} \sum_{z=0}^{\infty}\sum_{w=0}^{\infty}\sum_{u=0}^{\infty}(-1)^{z+w+u}\binom{\frac{-cb}{\alpha}}{z}\binom{c}{w}\binom{\frac{wb}{\alpha}}{u}\int_0^1 x^{a-1+ah+\frac{abl}{\alpha}+aj+am+ap} x^{az+a} dx$$

$$= W_{k,h,n,l,j,m,f,p} \left(\frac{\lambda}{\gamma}\right)^{2c} \sum_{z=0}^{\infty}\sum_{w=0}^{\infty}\sum_{u=0}^{\infty}(-1)^{z+w+u}\binom{\frac{-cb}{\alpha}}{z}\binom{c}{w}\binom{\frac{wb}{\alpha}}{u} \left.\frac{x^{a+ah+\frac{abl}{\alpha}+aj+am+ap+az}}{a+ah+\frac{abl}{\alpha}+aj+am+ap+az+au}\right|_0^1$$

$$E\left[\left(\frac{\lambda\left[(1-x^{a})^{\frac{-b}{\alpha}}-1\right]}{\gamma}\right)^{c}\right]$$

$$= W_{k,h,n,l,j,m,f,p} \left(\frac{\lambda}{\gamma}\right)^{2c} \sum_{z=0}^{\infty}\sum_{w=0}^{\infty}\sum_{u=0}^{\infty}(-1)^{z+w+u}\binom{\frac{-cb}{\alpha}}{z}\binom{c}{w}\binom{\frac{wb}{\alpha}}{u}\frac{1}{a\left(1+h+\frac{bl}{\alpha}+j+m+p+z+u\right)}$$

Put the results of the four parts together:

$$\frac{(a-1)\frac{abc}{\alpha}\left(\frac{\lambda}{\gamma}\right)^{c} W_{k,h,n,l,j,m,f,p}}{a\left(1+h+\frac{bl}{\alpha}+j+m+p\right)^{2}}$$

$$+\frac{\left(1+\frac{b}{\alpha}\right)\frac{abc}{\alpha}\left(\frac{\lambda}{\gamma}\right)^{c} W_{k,h,n,l,j,m,f,p}}{\left(1+h+\frac{bl}{\alpha}+j+m+p\right)}\sum_{s=0}^{\infty}\frac{1}{a\left(2+h+\frac{bl}{\alpha}+j+m+p+s\right)}$$

$$+\frac{b(c-1)\frac{abc}{\alpha}\left(\frac{\lambda}{\gamma}\right)^{c} W_{k,h,n,l,j,m,f,p}}{a\left(1+h+\frac{bl}{\alpha}+j+m+p\right)}\sum_{v=0}^{\infty}\sum_{r=0}^{\infty}\sum_{d=0}^{\infty}(-1)^{d}\binom{\frac{br}{\alpha}}{d}\frac{1}{a\left(2+h+\frac{bl}{\alpha}+j+m+p+v+d\right)}$$

$$+ W_{k,h,n,l,j,m,f,p} \frac{abc}{\alpha}\left(\frac{\lambda}{\gamma}\right)^{2c} \sum_{z=0}^{\infty}\sum_{w=0}^{\infty}\sum_{u=0}^{\infty}(-1)^{z+w+u}\binom{\frac{-cb}{\alpha}}{z}\binom{c}{w}\binom{\frac{wb}{\alpha}}{u}.\frac{1}{a\left(1+h+\frac{bl}{\alpha}+j+m+p+z+u\right)}$$

Let us call: $\varphi = 1+h+\frac{bl}{\alpha}+j+m+p$ ,

$$W_1 = \frac{abc}{\alpha}\left(\frac{\lambda}{\gamma}\right)^c W_{k,h,n,l,j,m,f,p} \quad , \quad W_2 = \sum_{v=0}^{\infty}\sum_{r=0}^{\infty}\sum_{d=0}^{\infty}(-1)^d \binom{\frac{br}{\alpha}}{d} = W_{v,r,d}$$

$$W_3 = \sum_{z=0}^{\infty}\sum_{w=0}^{\infty}\sum_{u=0}^{\infty}(-1)^{z+w+u}\binom{\frac{-cb}{\alpha}}{z}\binom{c}{w}\binom{\frac{wb}{\alpha}}{u} = W_{z,w,u}$$

So we can rewrite the four above quantities in compact form

$$\frac{(a-1)W_1}{a(\varphi)^2} + \frac{\left(1+\frac{b}{\alpha}\right)W_1}{(\varphi)} \sum_{s=0}^{\infty}\frac{1}{a(1+\varphi+s)} + \frac{b(c-1)W_1\, W_2}{a(\varphi)}\frac{1}{a(1+\varphi+v+d)}$$

$$+ W_1 \left(\frac{\lambda}{\gamma}\right)^c W_3 \frac{1}{a(\varphi+z+u)}$$

$$\therefore E\{-\log(g[x])\} = -\ln a - \ln b - \ln c - c\, \ln\lambda + \ln\alpha + c\ln\gamma + E\{-(a-1)\ln x\}$$

$$+\left(1+\frac{b}{\alpha}\right) E\{\ln(1-x^a)\} + E\left\{-(c-1)\ln\left[(1-x^a)^{\frac{-b}{\alpha}} - 1\right]\right\}$$

$$+ E\left\{\left(\frac{\lambda}{\gamma}\right)^c \left((1-x^a)^{\frac{-b}{\alpha}} - 1\right)^c\right\}$$

$$\therefore E\{-\log(g[x])\} = -\ln a - \ln b - \ln c - c\, \ln\lambda + \ln\alpha + c\ln\gamma + \frac{(a-1)W_1}{a(\varphi)^2}$$

$$+ \frac{\left(1+\frac{b}{\alpha}\right)W_1}{(\varphi)} \sum_{s=0}^{\infty}\frac{1}{a(1+\varphi+s)} + \frac{b(c-1)W_1\, W_2}{a(\varphi)}\frac{1}{a(1+\varphi+v+d)}$$

$$+ W_1 \left(\frac{\lambda}{\gamma}\right)^c W_3 \frac{1}{a(\varphi+z+u)}$$

## 4. Figures of The PDF, CDF Quantile and Hazard Function

In this section, the author depicts different figures (1-12) illustrating the behavior of the distribution. The PDF can exhibit bathtub, J shape and reversed J shaped, increasing and decreasing, N shape and other different shapes.

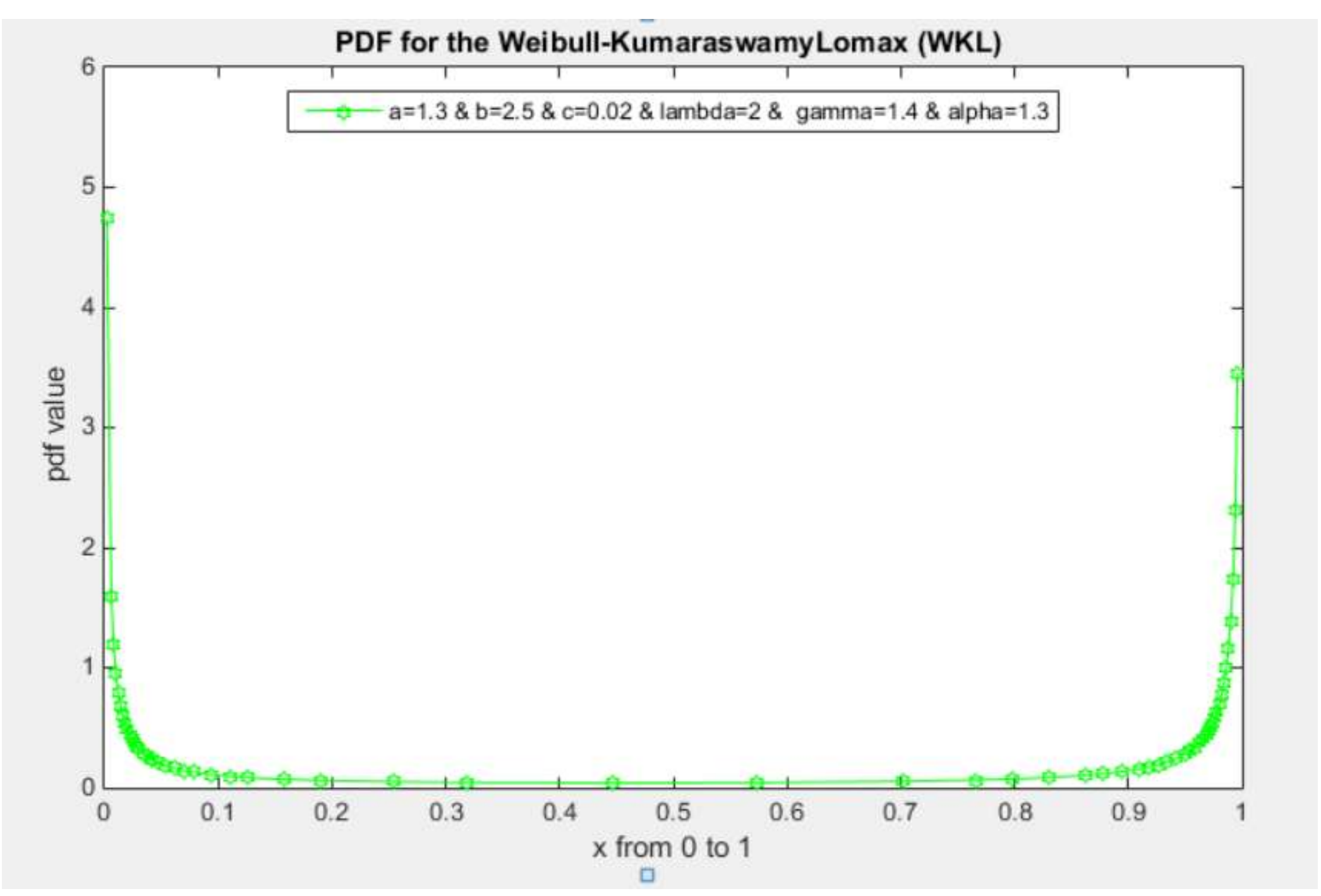

Fig. 1 shows the bathtub shape of the PDF with a=1.3, b=2.5, c=0.02, lambda=2, gamma=1.4, alpha=1.3

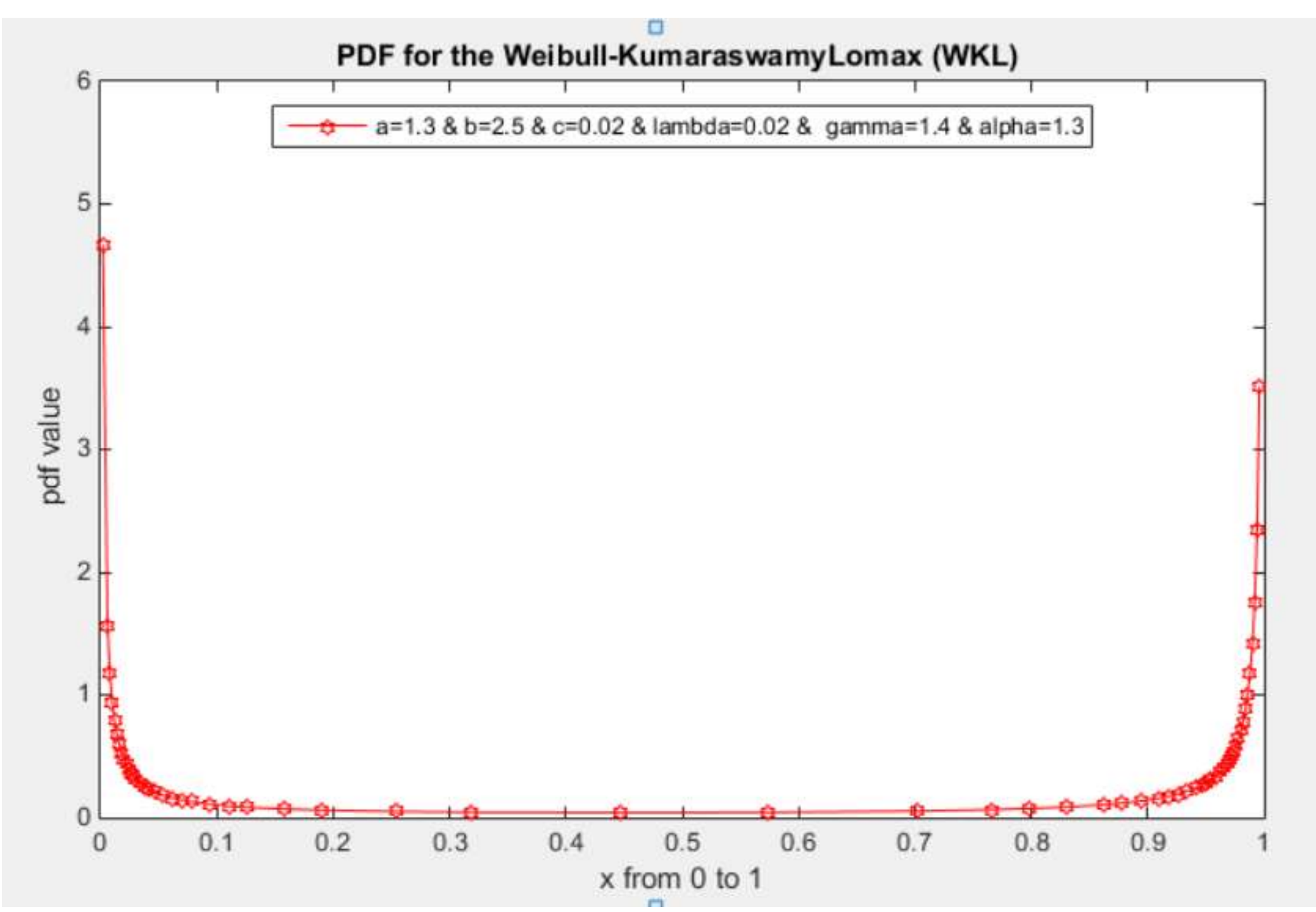

Fig. 2 shows the bathtub shape of the PDF with a=1.3, b=2.5, c=0.02, lambda=0.02, gamma=1.4, alpha=1.3. Changing the lambda value from 2 in fig.1 to 0.02 in fig.2 while keeping other parameters with the same values did not change the bathtub shape of the PDF.

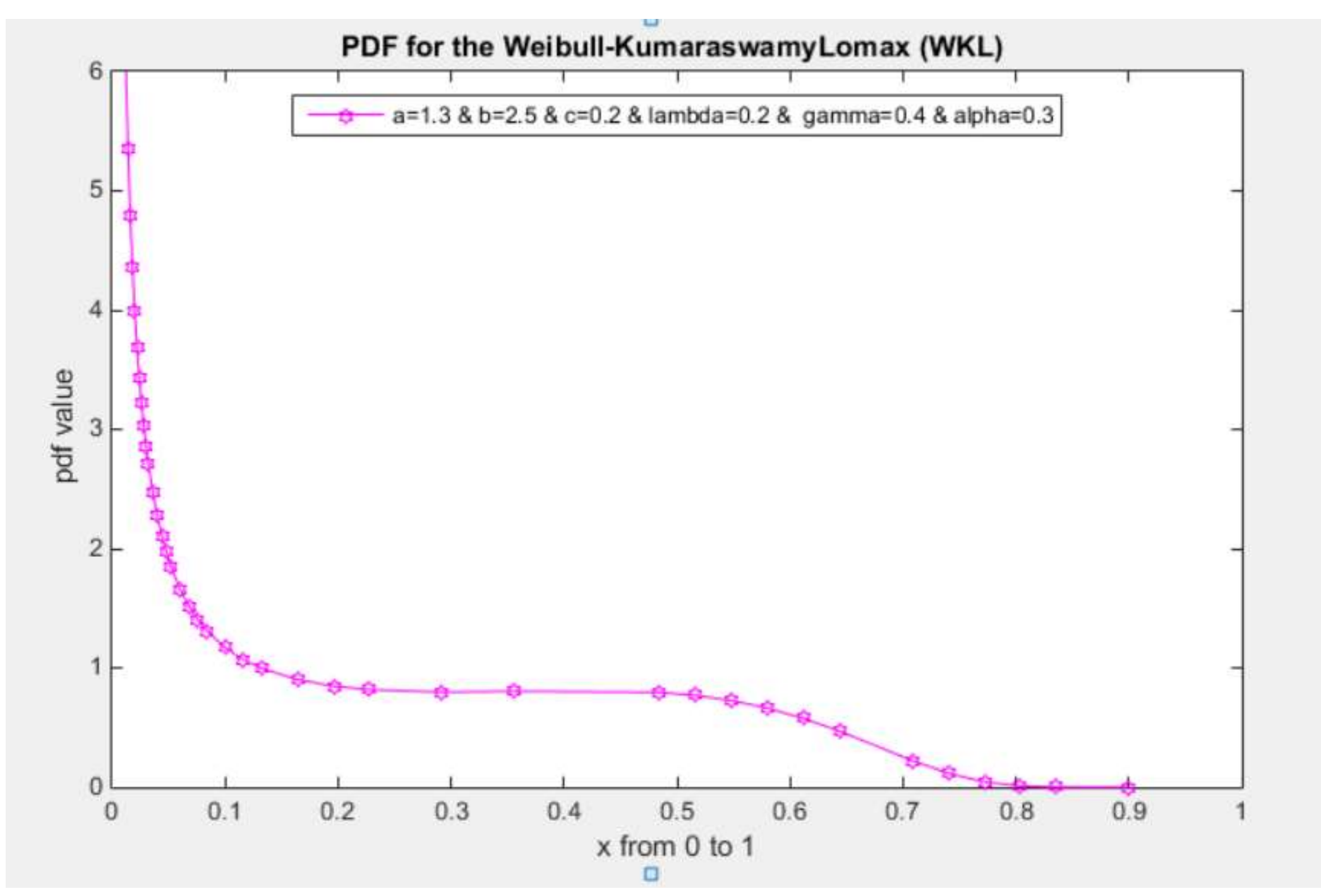

Fig. 3 shows the decreasing shape of the PDF with a=1.3, b=2.5, c=0.2, lambda=0.2, gamma=0.4, alpha=0.3.

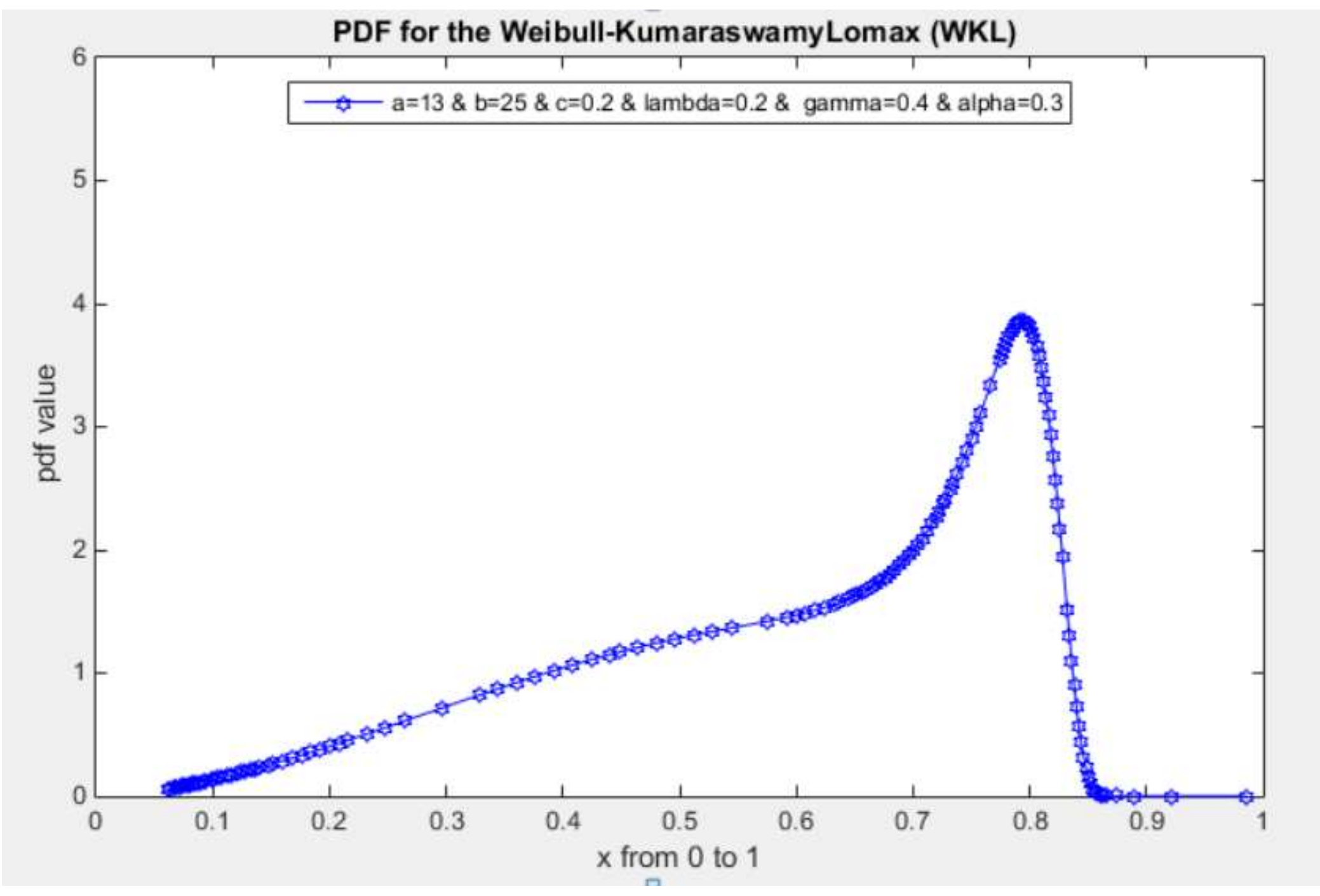

Fig. 4 shows the increasing then the decreasing shape of the PDF with a=13, b=25, c=0.2, lambda=0.2, gamma=0.4, alpha=0.3.

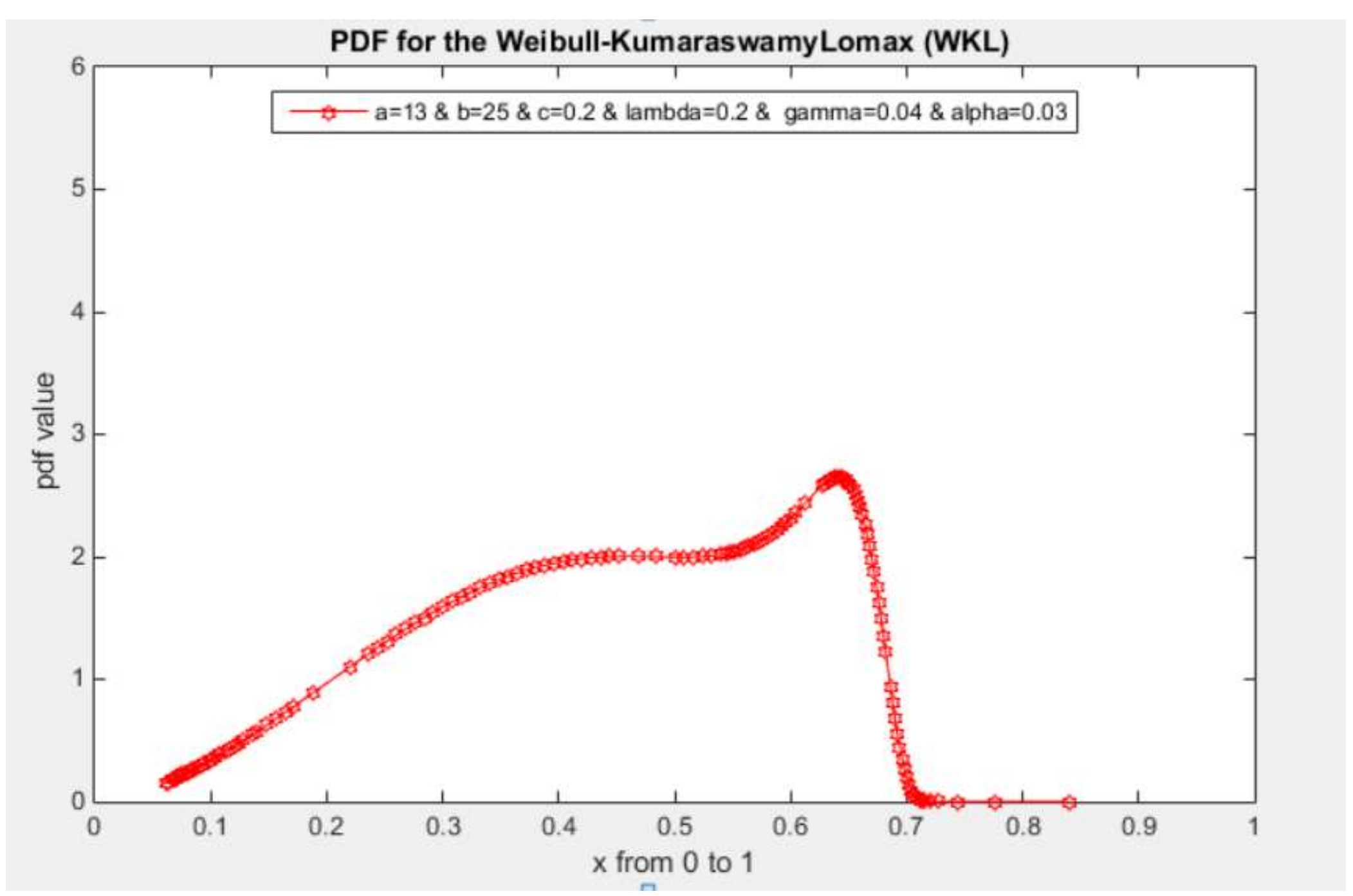

Fig. 5 shows the increasing then the decreasing shape of the PDF with a=13, b=25, c=0.2, lambda=0.2, gamma=0.04, alpha=0.03. The shape differs than in fig.4 when changing the values of the gamma and alpha.

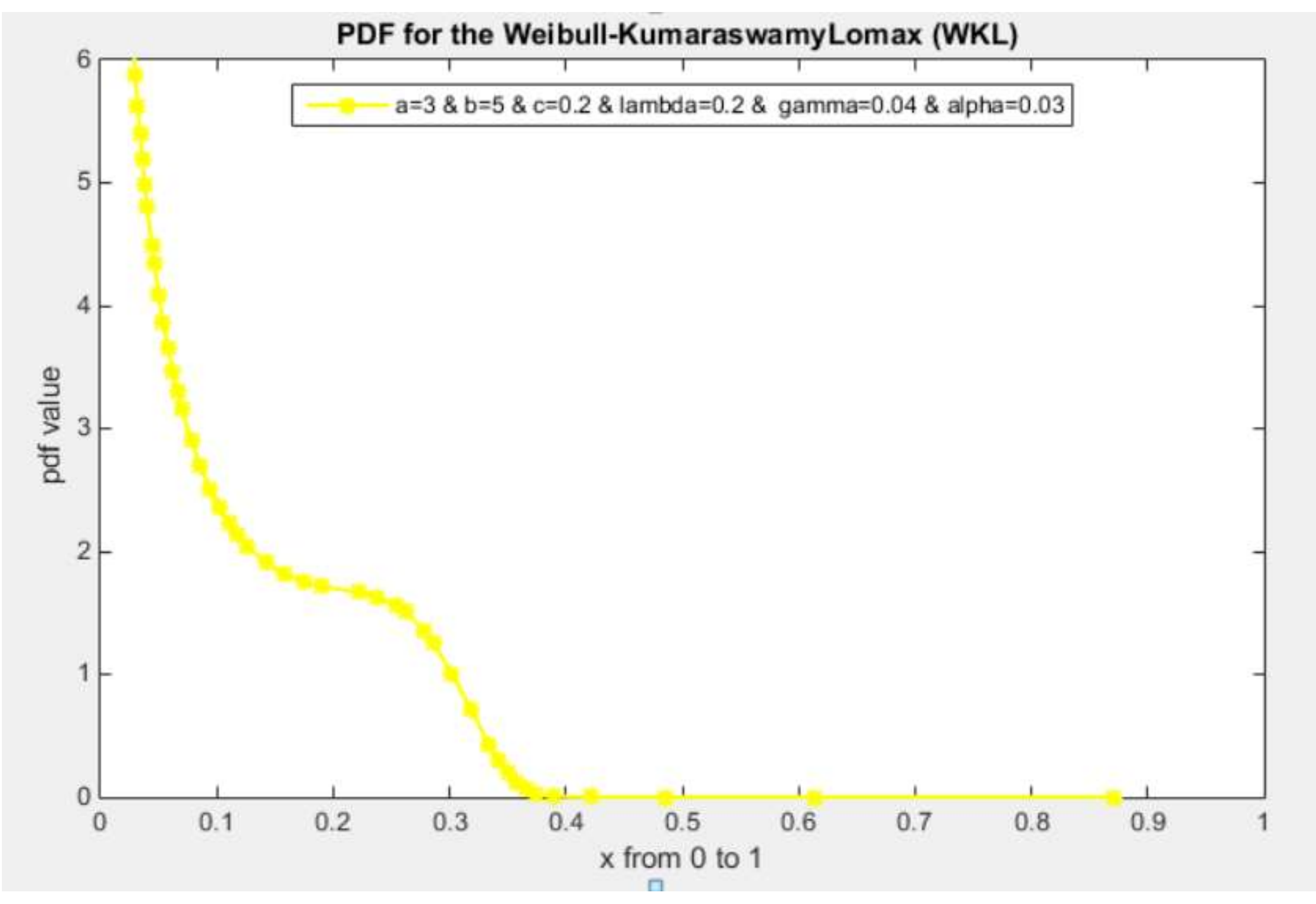

Fig. 6 shows the decreasing shape of the PDF and the associated hump before decreasing to the near zero values with a=3, b=5, c=0.2, lambda=0.2, gamma=0.04, alpha=0.03. The shape differs than in fig.5 when decreasing the values of a and b in fig 6.

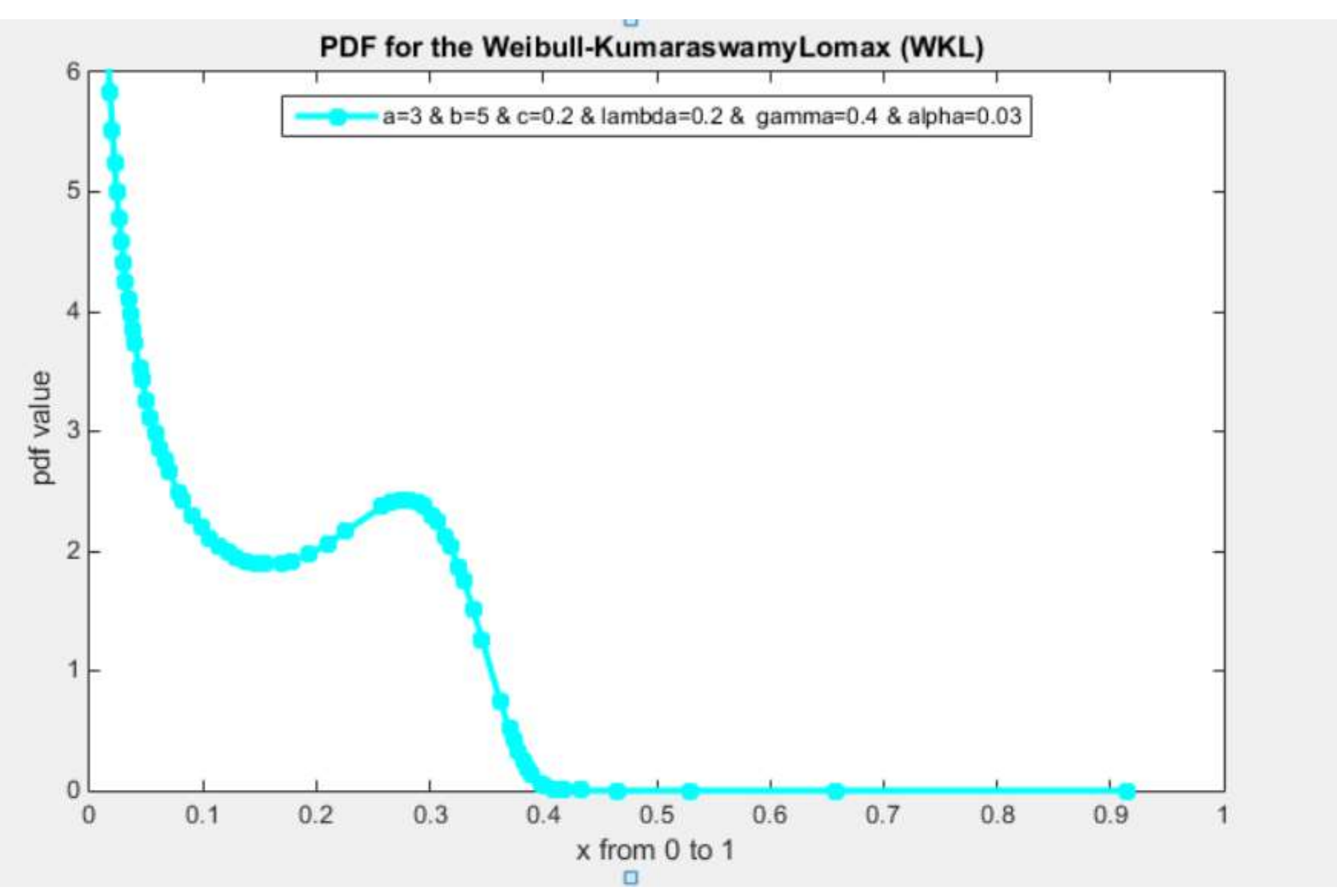

Fig. 7 shows the decreasing shape of the PDF and the associated hump before decreasing to near zero values with a=3, b=5, c=0.2, lambda=0.2, gamma=0.4, alpha=0.03. The shape differs than in fig.6 when increasing the values of gamma in fig 7. from 0.04 in fig. 6 to 0.4 in fig. 7. The hump is mainly positioned more at the left tail of the distribution.

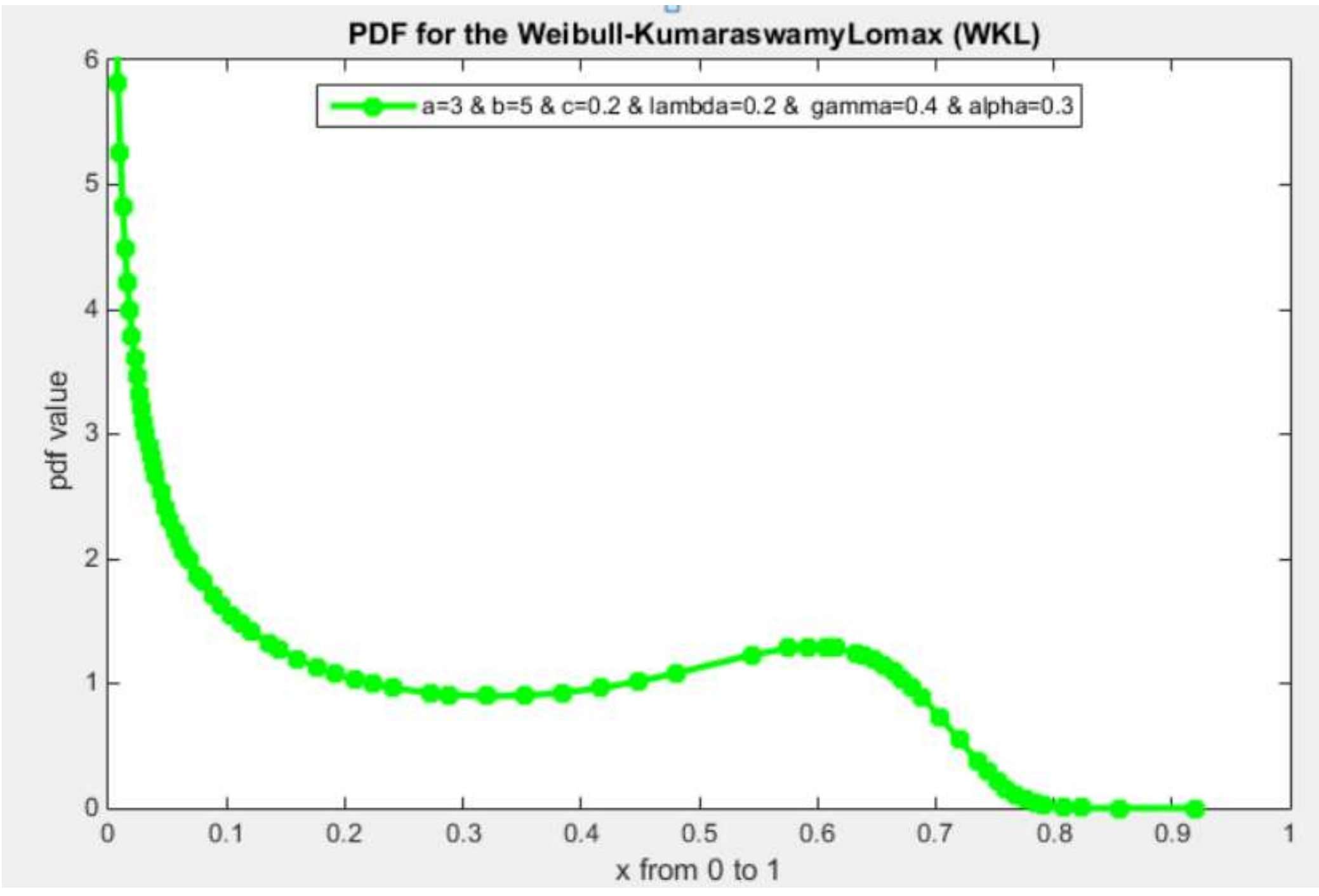

Fig. 8 shows the decreasing shape of the PDF and the associated hump before decreasing to near zero values with a=3, b=5, c=0.2, lambda=0.2, gamma=0.4, alpha=0.3. The shape differs than in fig.7 when increasing the values of alpha in fig 7. from 0.03 in fig. 7 to 0.3 in fig. 8. Also the hump changes its position to the right tail of the distribution.

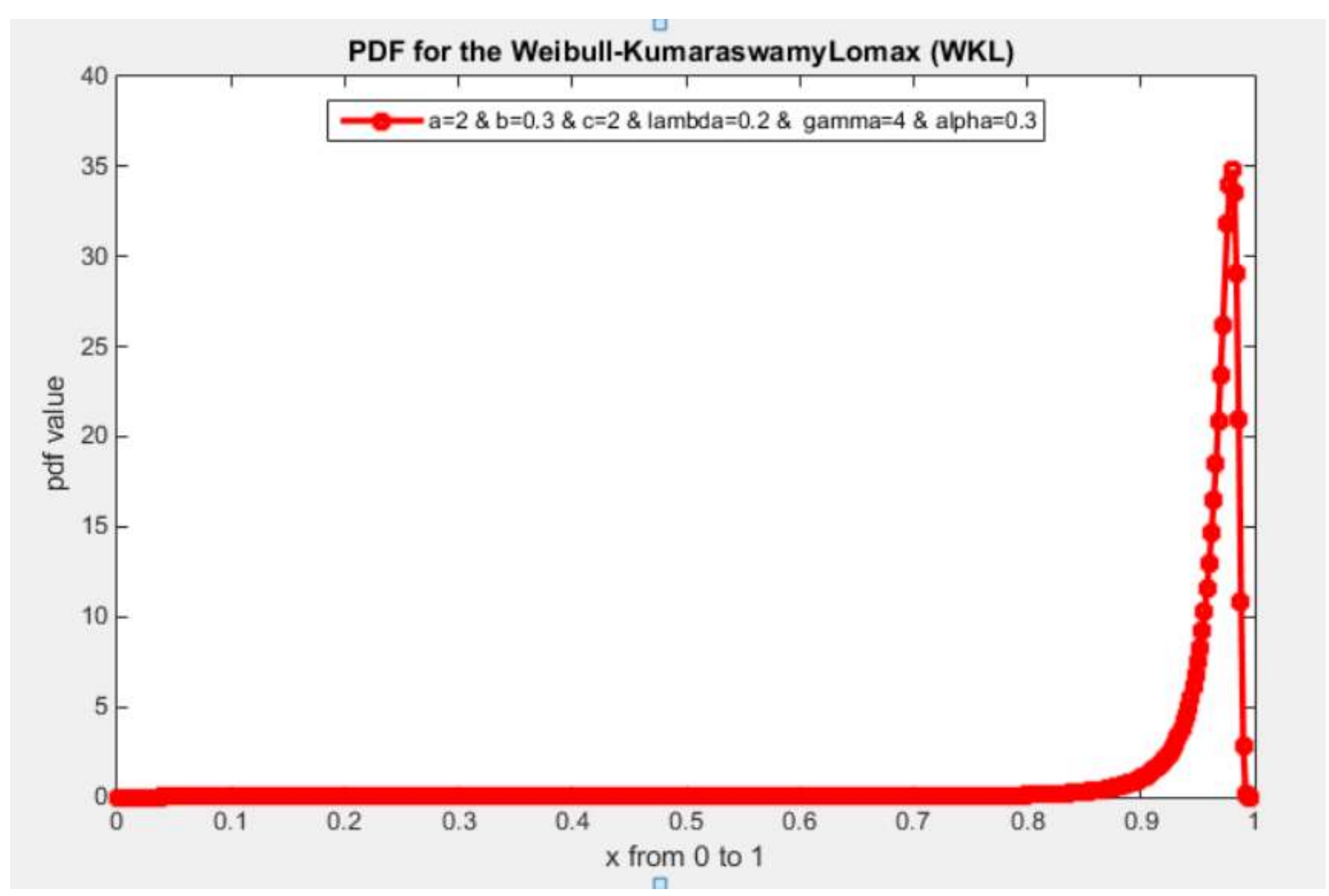

Fig. 9 shows the shape of the PDF mainly at the right tail upside and then downside with a=2, b=0.3, c=2, lambda=0.2, gamma=4, alpha=0.3.

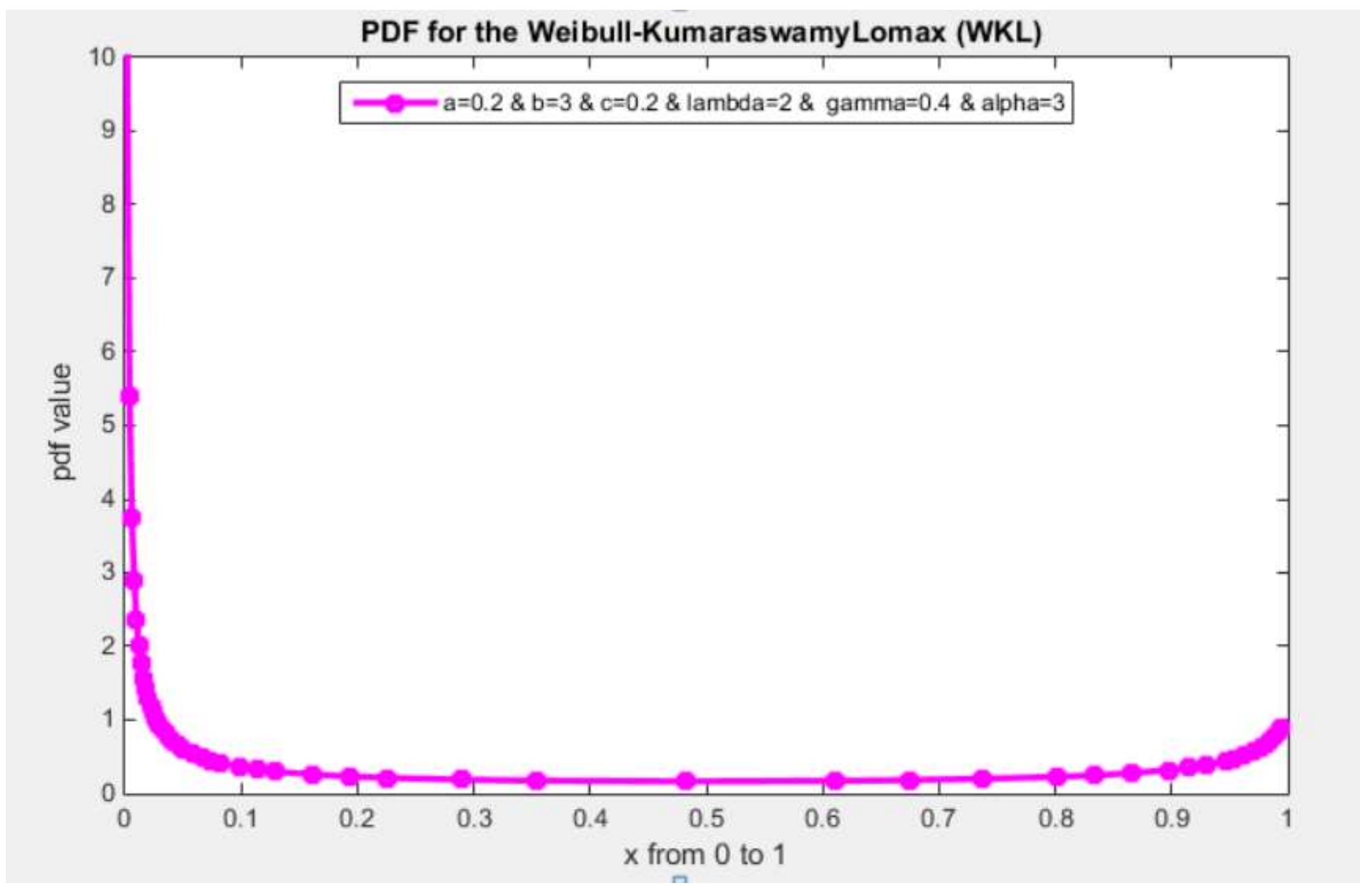

Fig. 10 shows the reversed J shaped of the PDF with a=0.2, b=3, c=0.2, lambda=2, gamma=0.4, alpha=3.

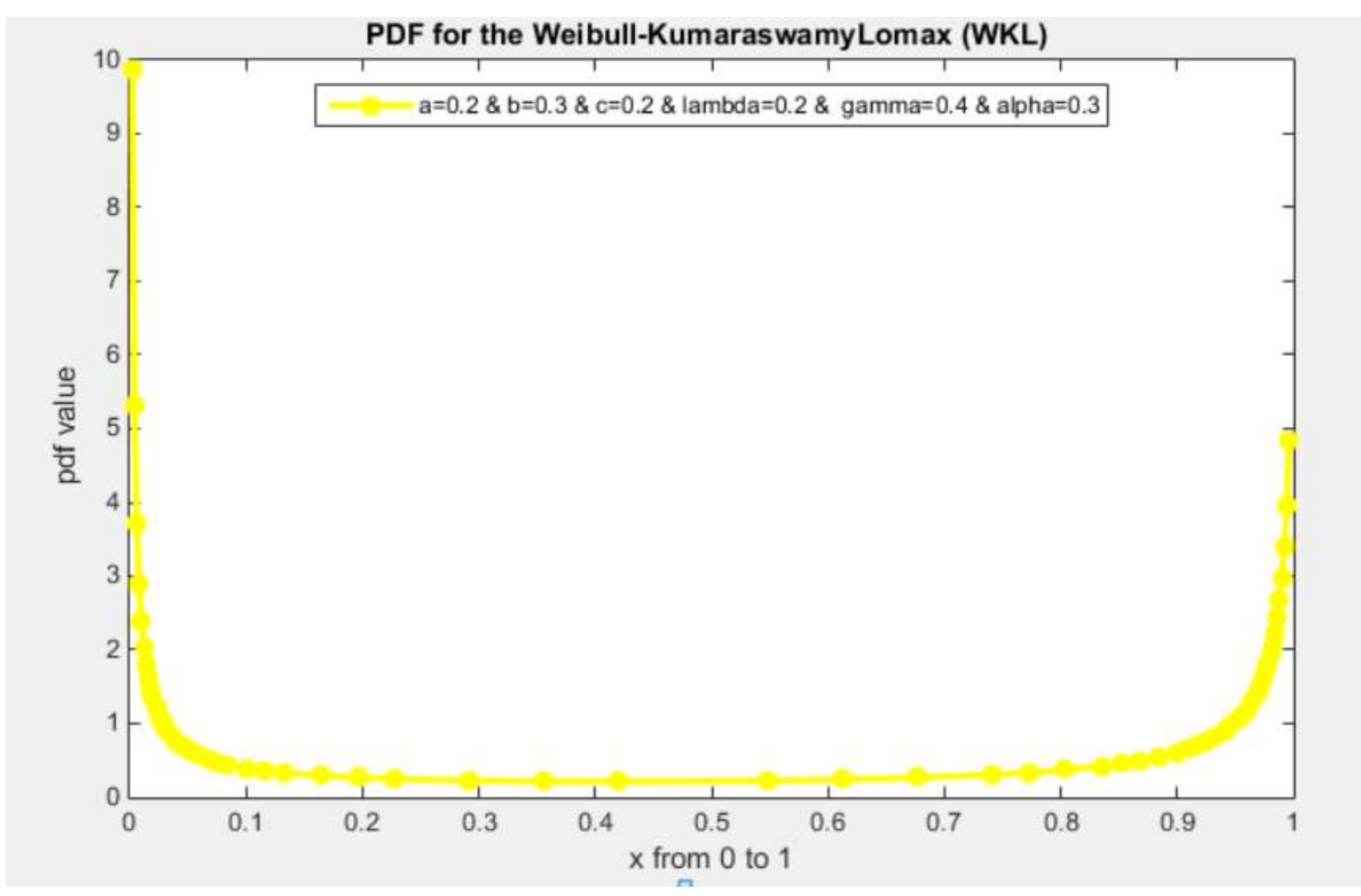

Fig. 11 shows the reversed J shape of the PDF with a=0.2, b=0.3, c=0.2, lambda=0.2, gamma=0.4, alpha=0.3.

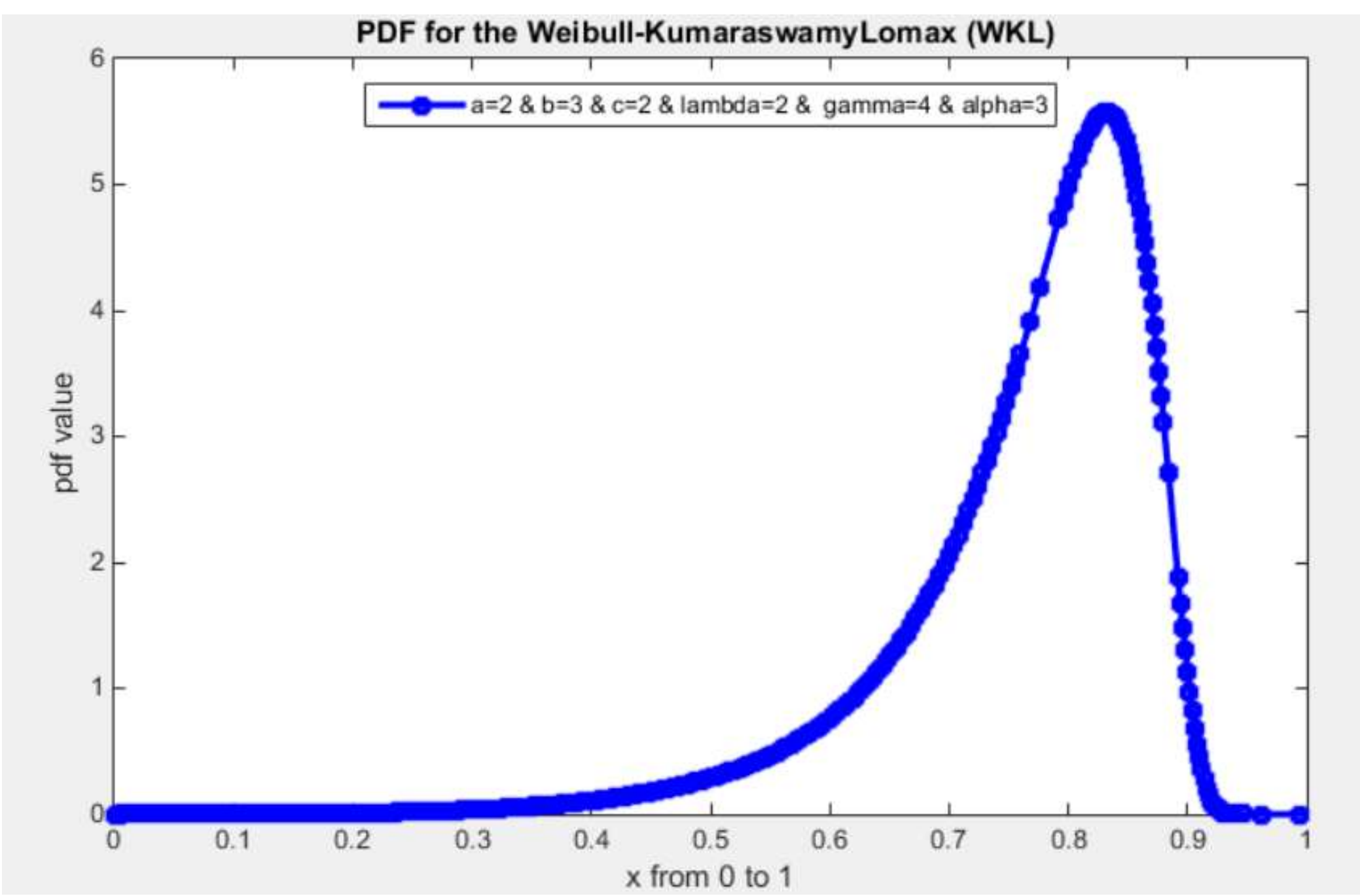

Fig. 12 shows the upside and downside of the PDF with a=2, b=3, c=2, lambda=2, gamma=4, alpha=3 mainly at the right tail of the distribution.

The Hazard function of the new distribution can mainly exhibit bathtub shape or increasing pattern as illustrated in the following figures (13-24) depending on the values of the parameters. In each figure the values of the parameters are illustrated at the top of the figure.

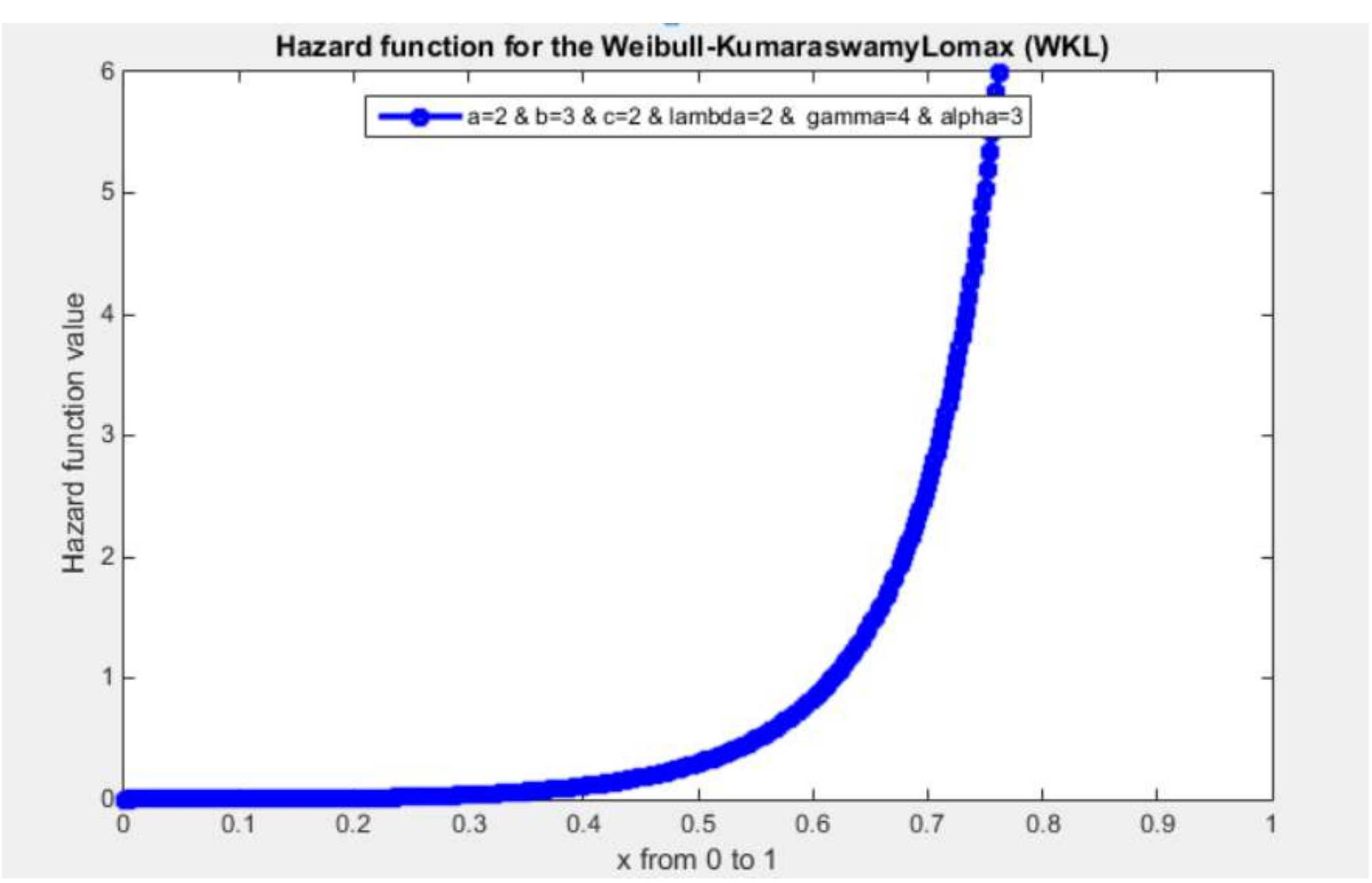

Fig. 13 shows increasing hazard rate function at the specified parameters values

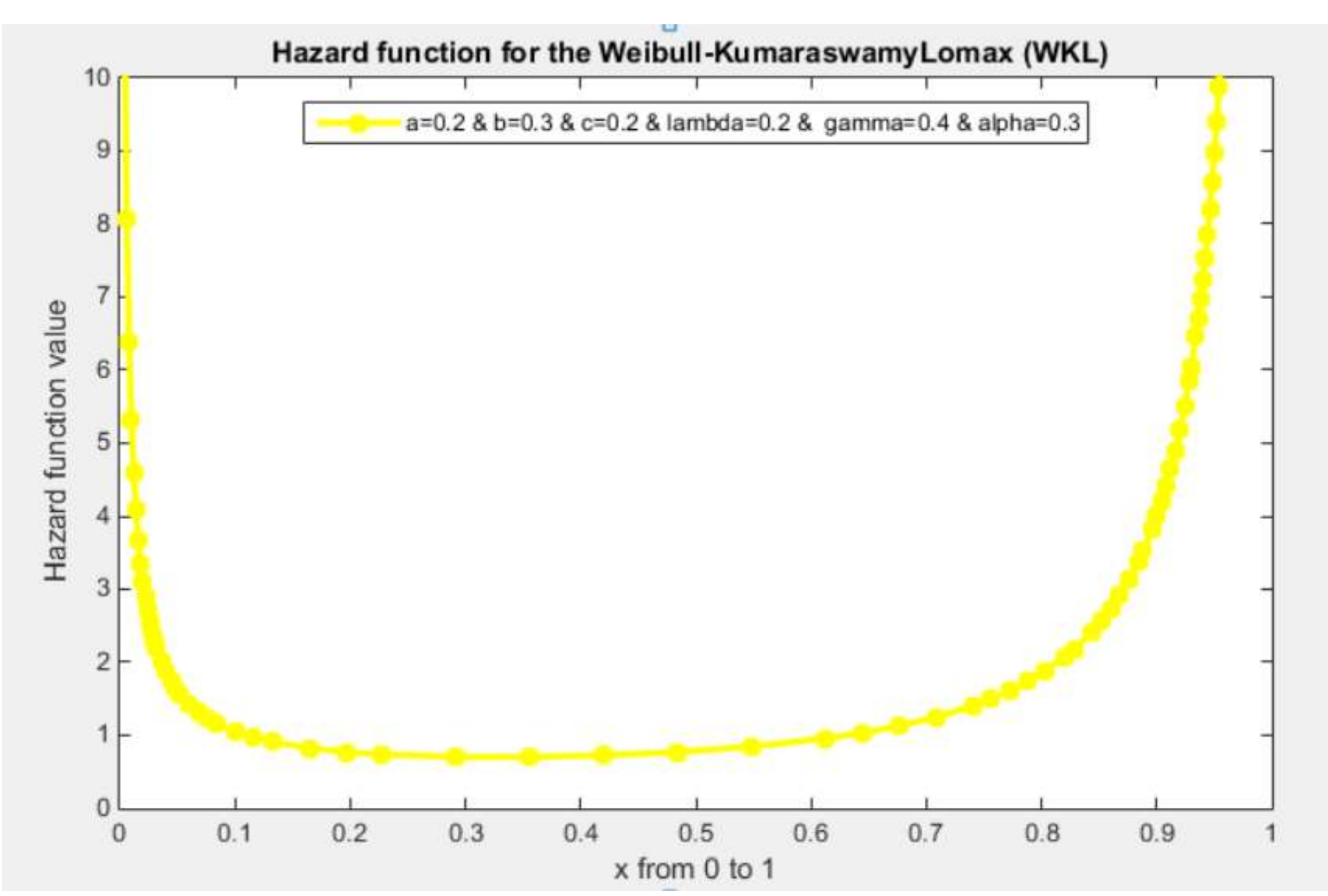

Fig. 14 shows bathtub hazard rate function at the specified parameters values.

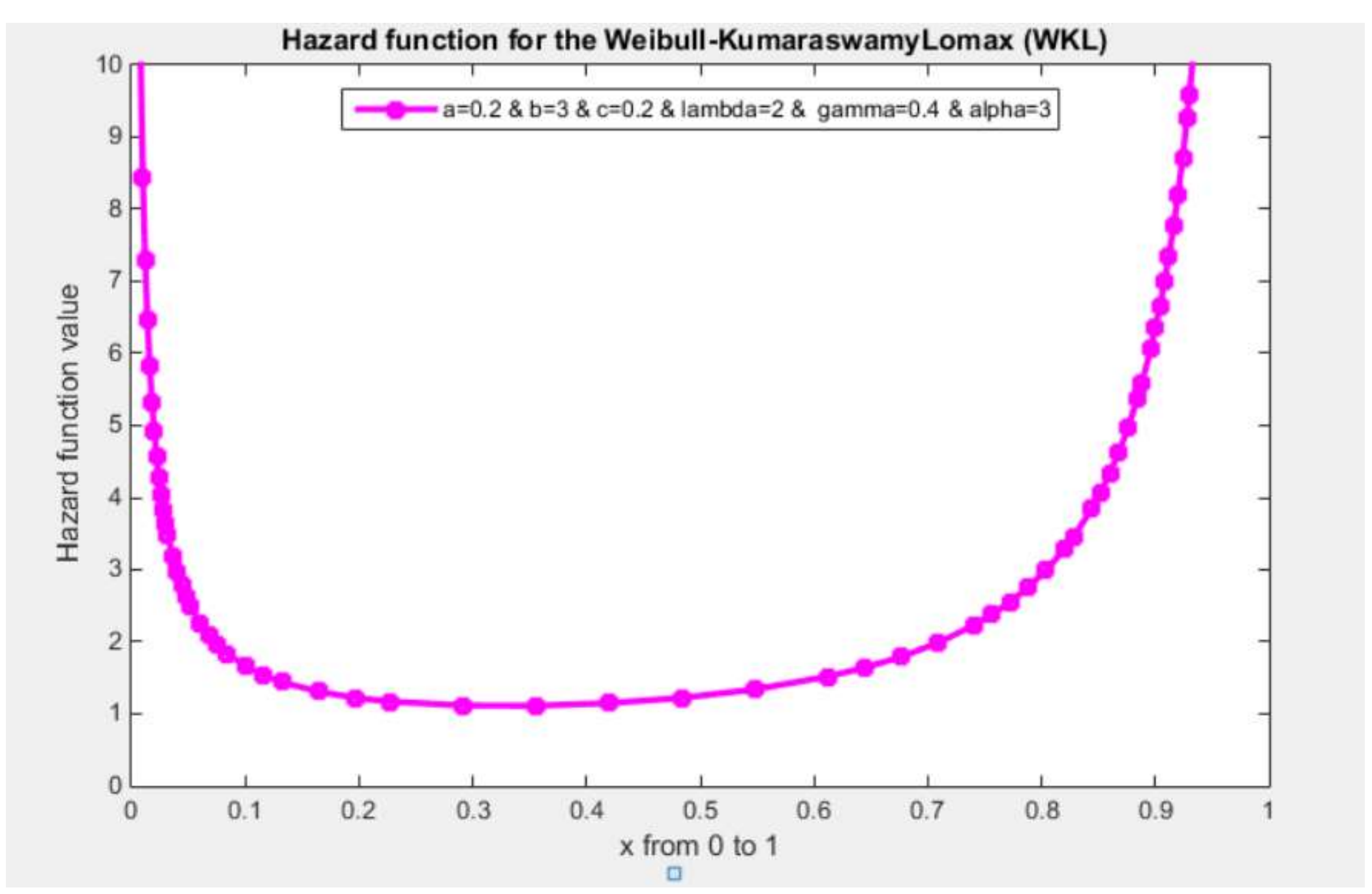

Fig. 15 shows bathtub hazard rate function at the specified parameters values

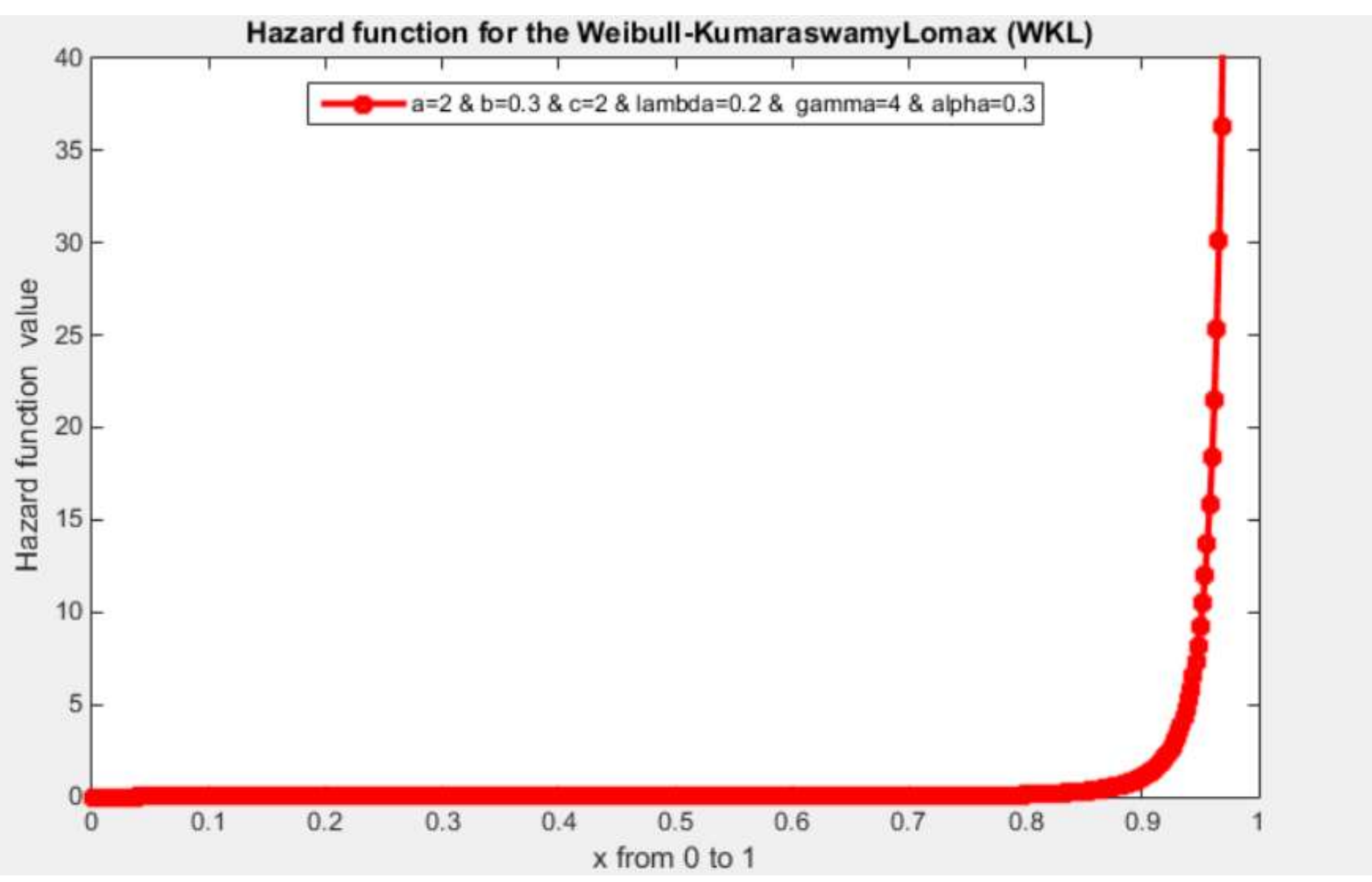

Fig. 16 shows increasing hazard rate function at the specified parameters values.

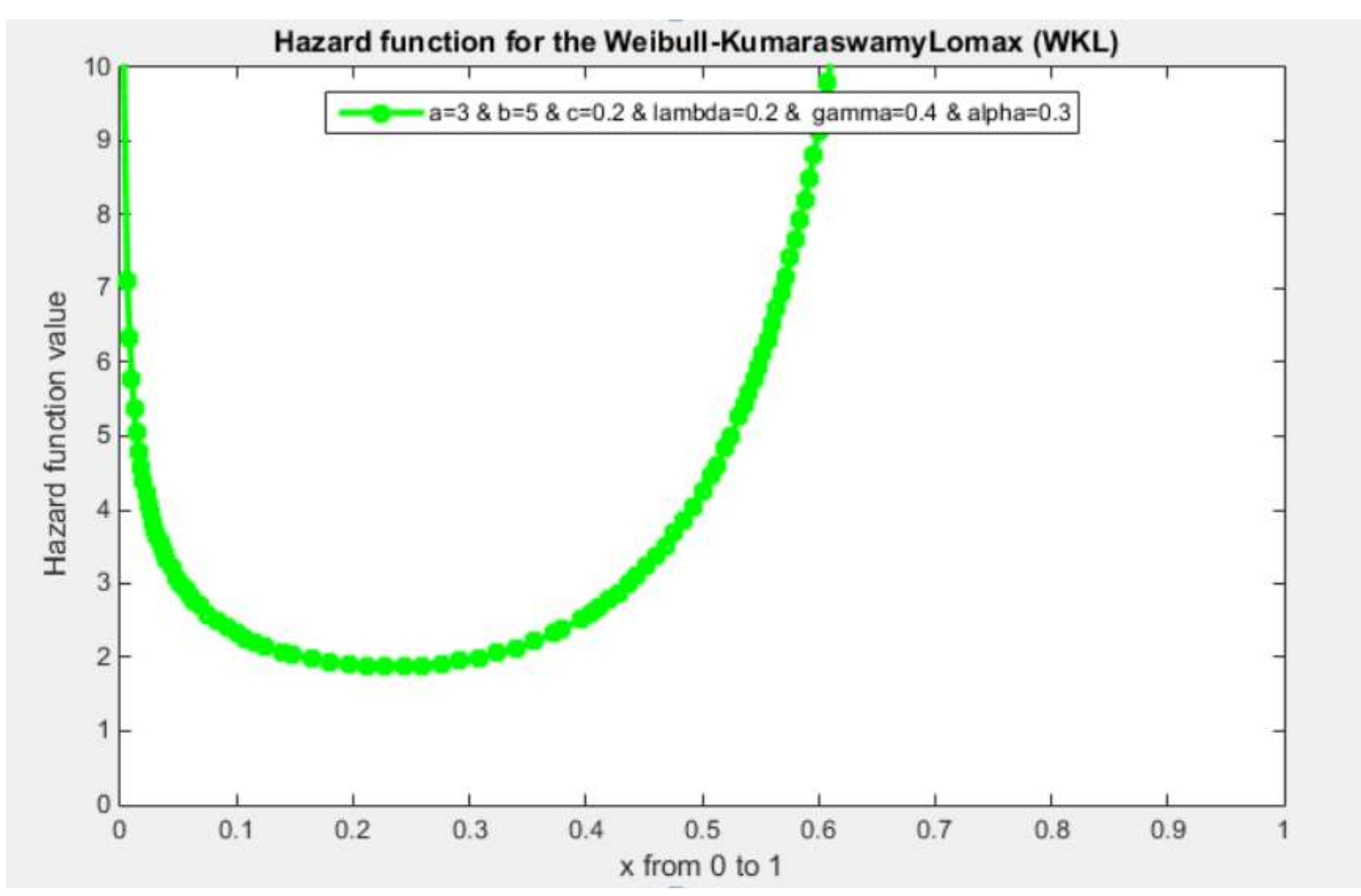

Fig. 17 shows bathtub hazard rate function at the specified parameters values.

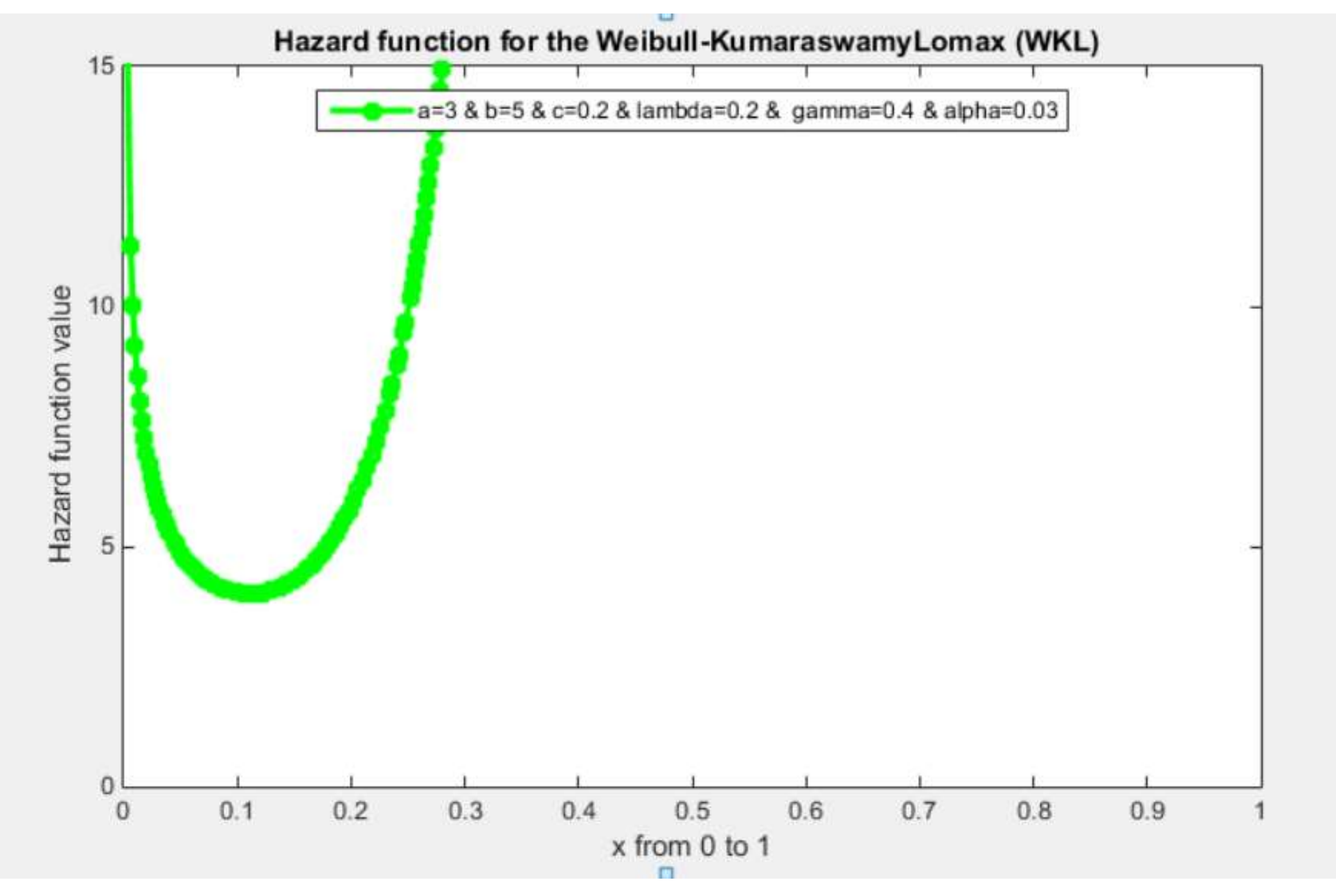

Fig. 18 shows bathtub hazard rate function at the specified parameters values.

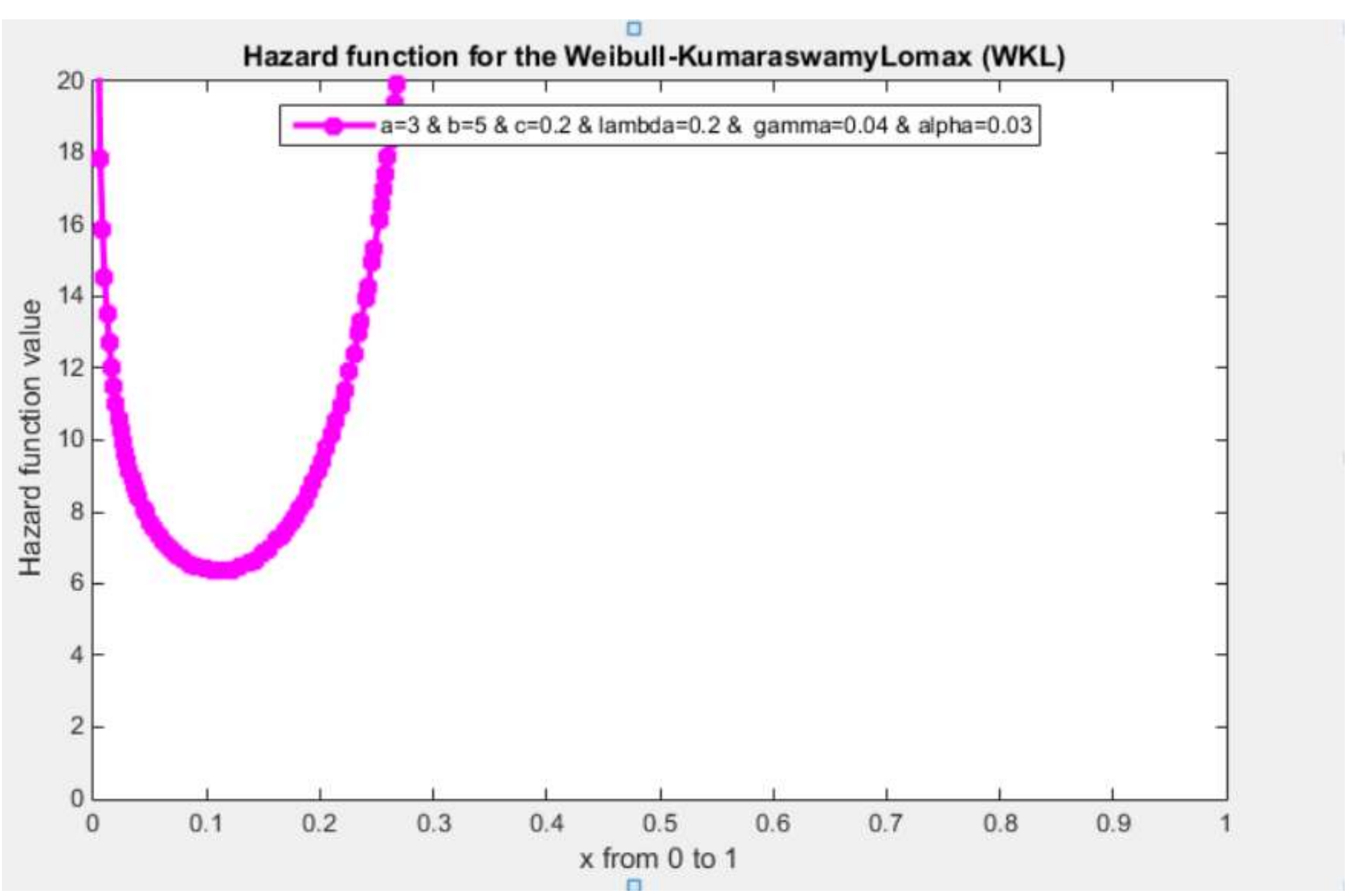

Fig. 19 shows bathtub hazard rate function at the specified parameters values.

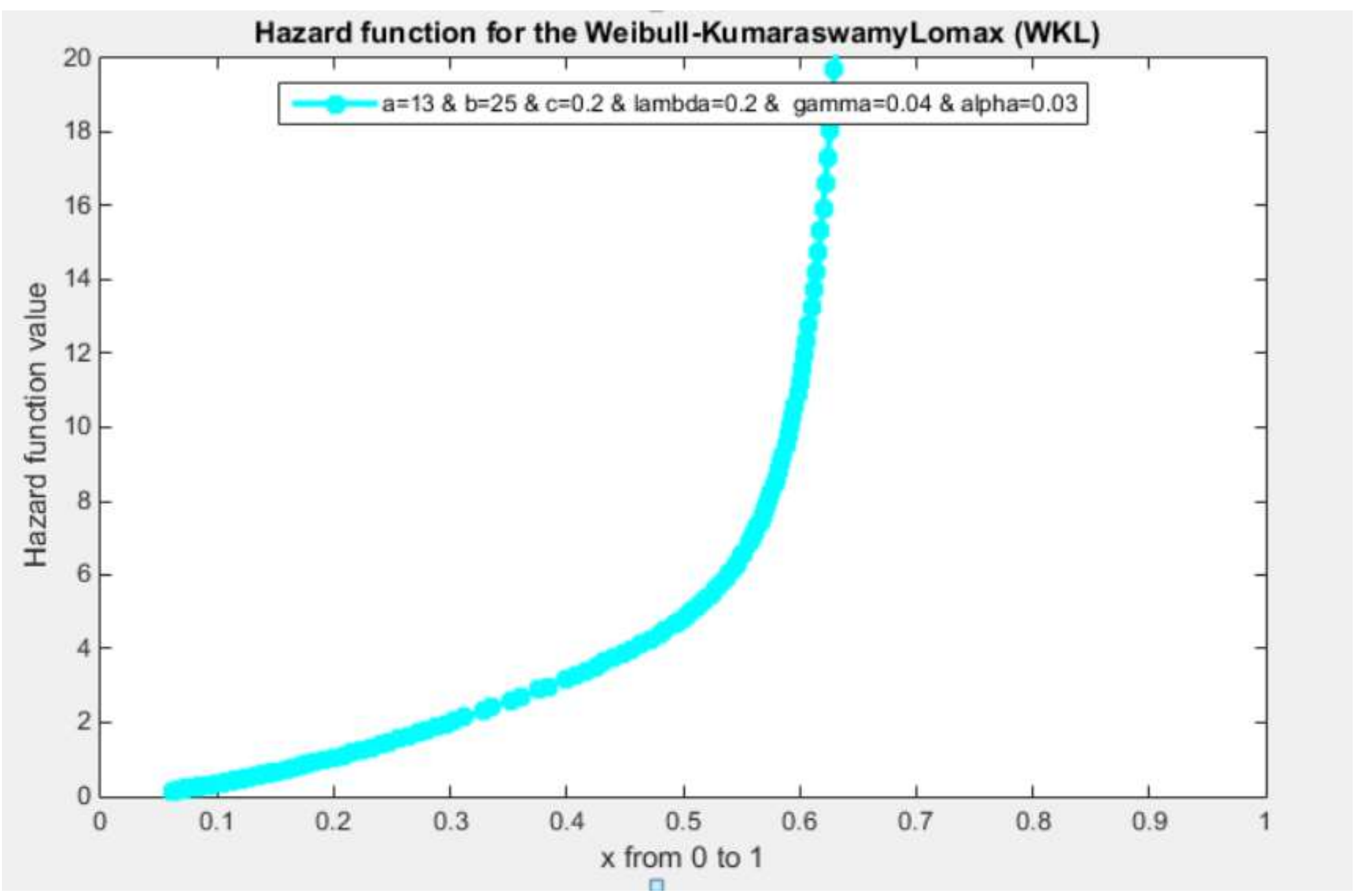

Fig. 20 shows increasing hazard rate function at the specified parameters values.

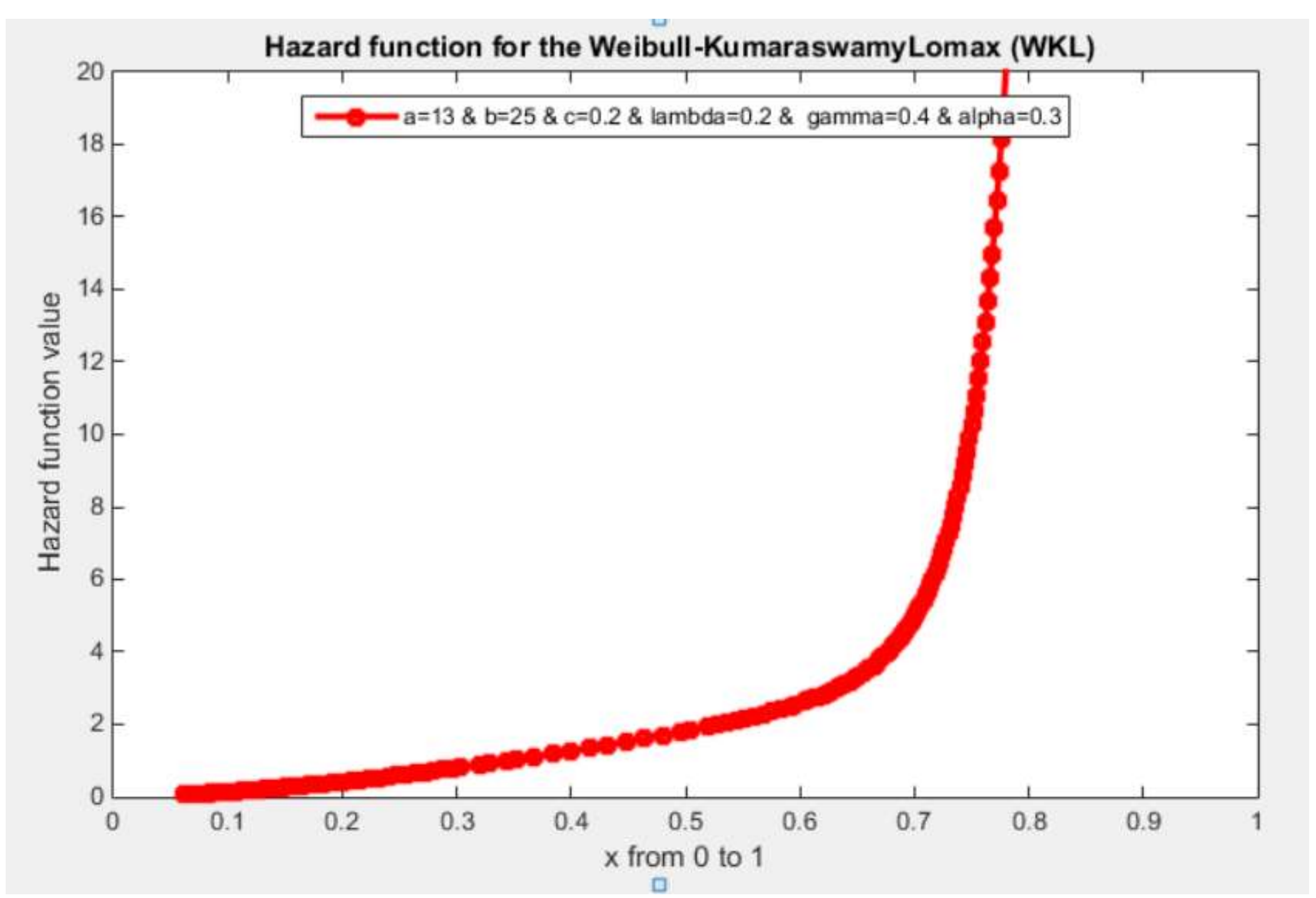

Fig. 21 shows increasing hazard rate function at the specified parameters values.

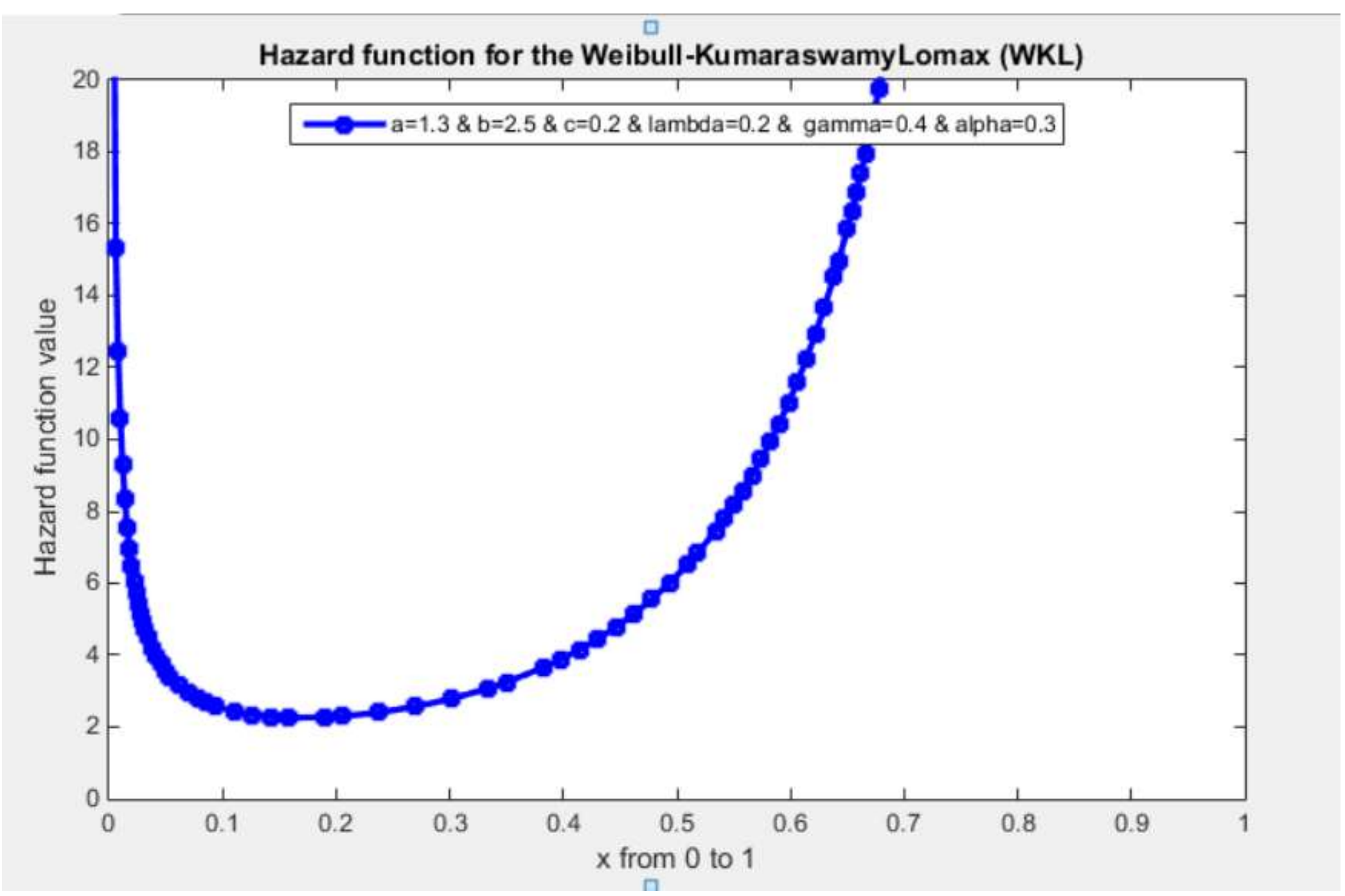

Fig. 22 shows bathtub hazard rate function at the specified parameters values.

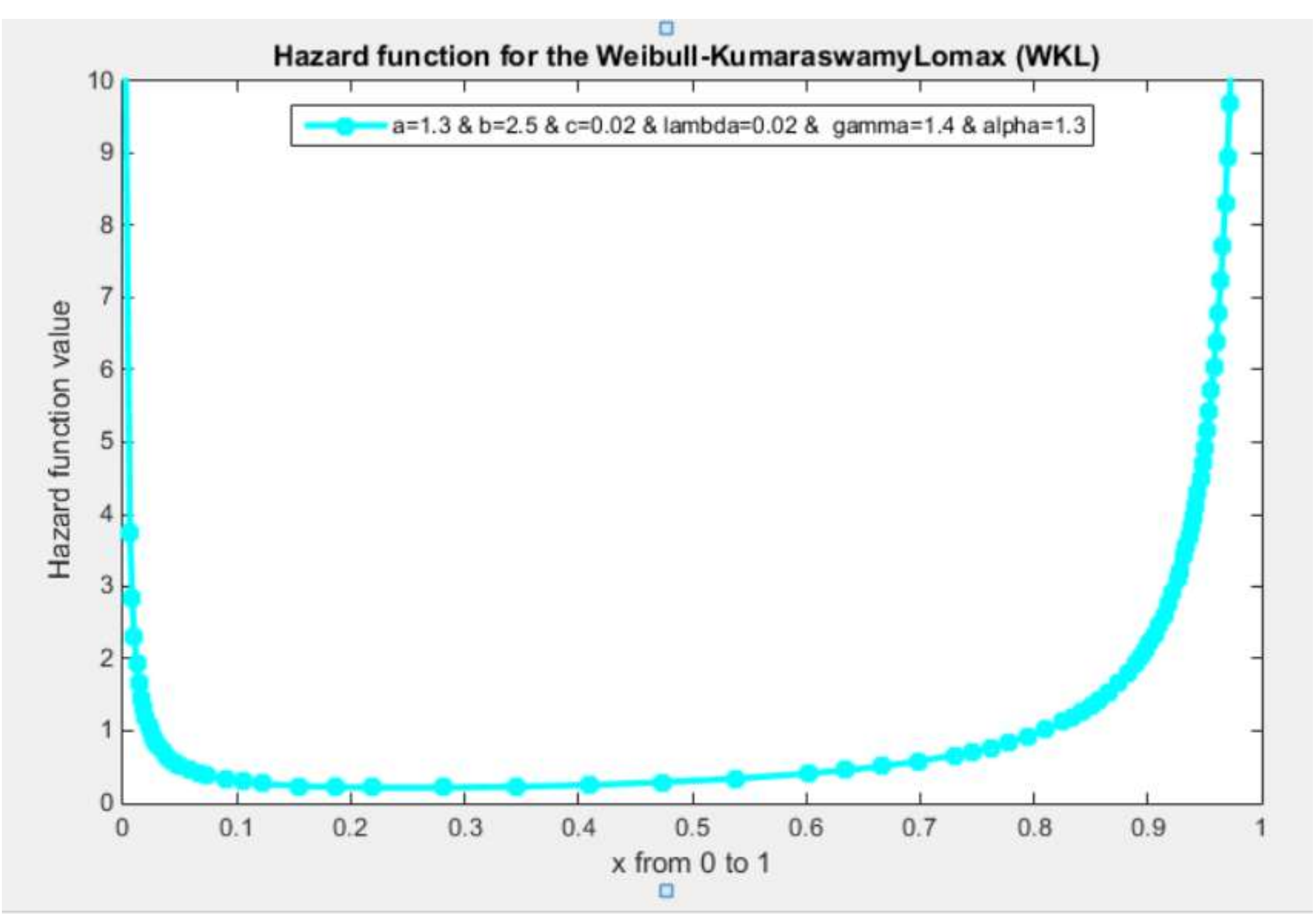

Fig. 23 shows bathtub hazard rate function at the specified parameters values.

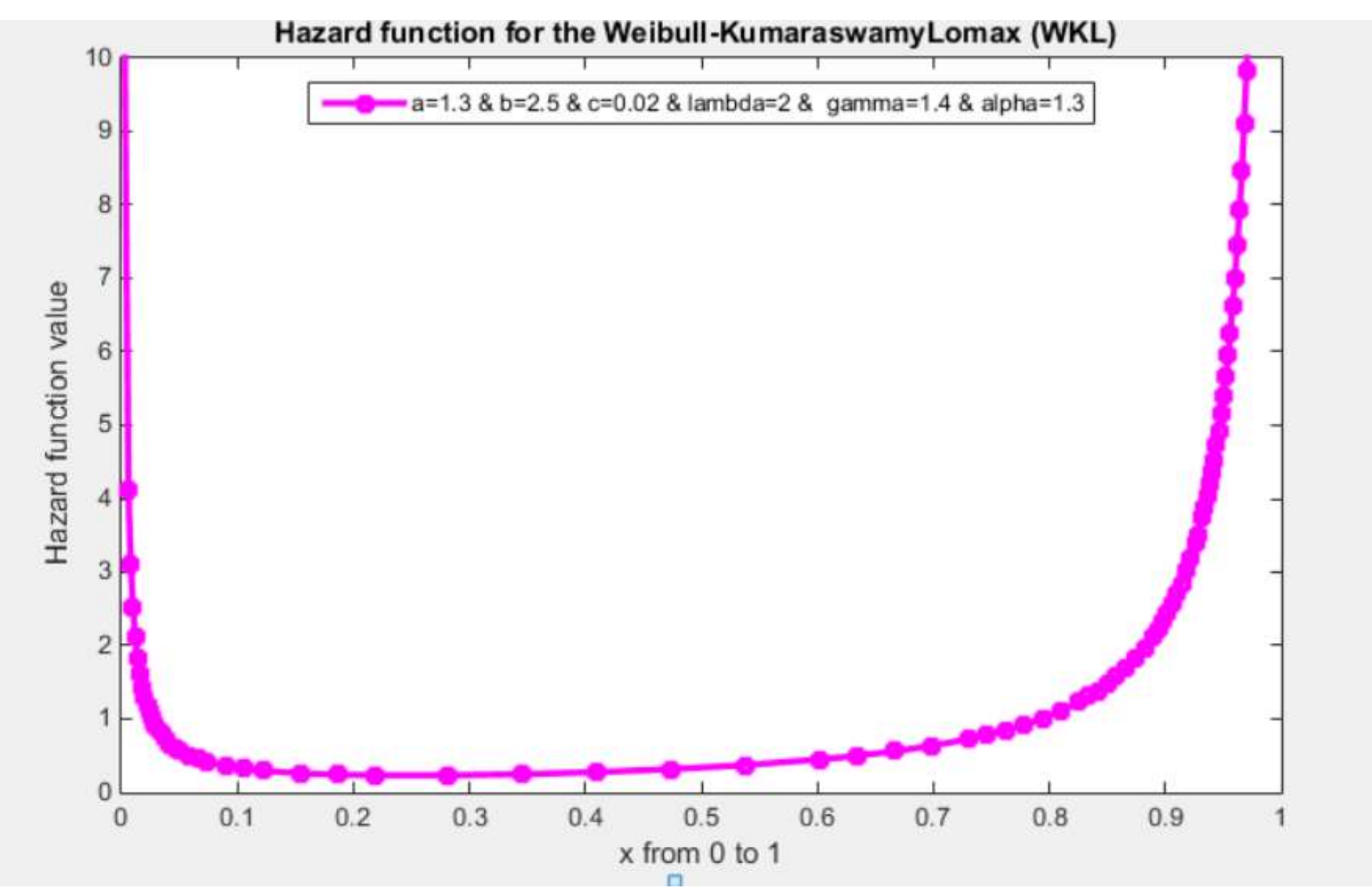

Fig. 24 shows the bathtub hazard rate function at the specified parameters values.

The CDF is shown in the following figures (25-27).

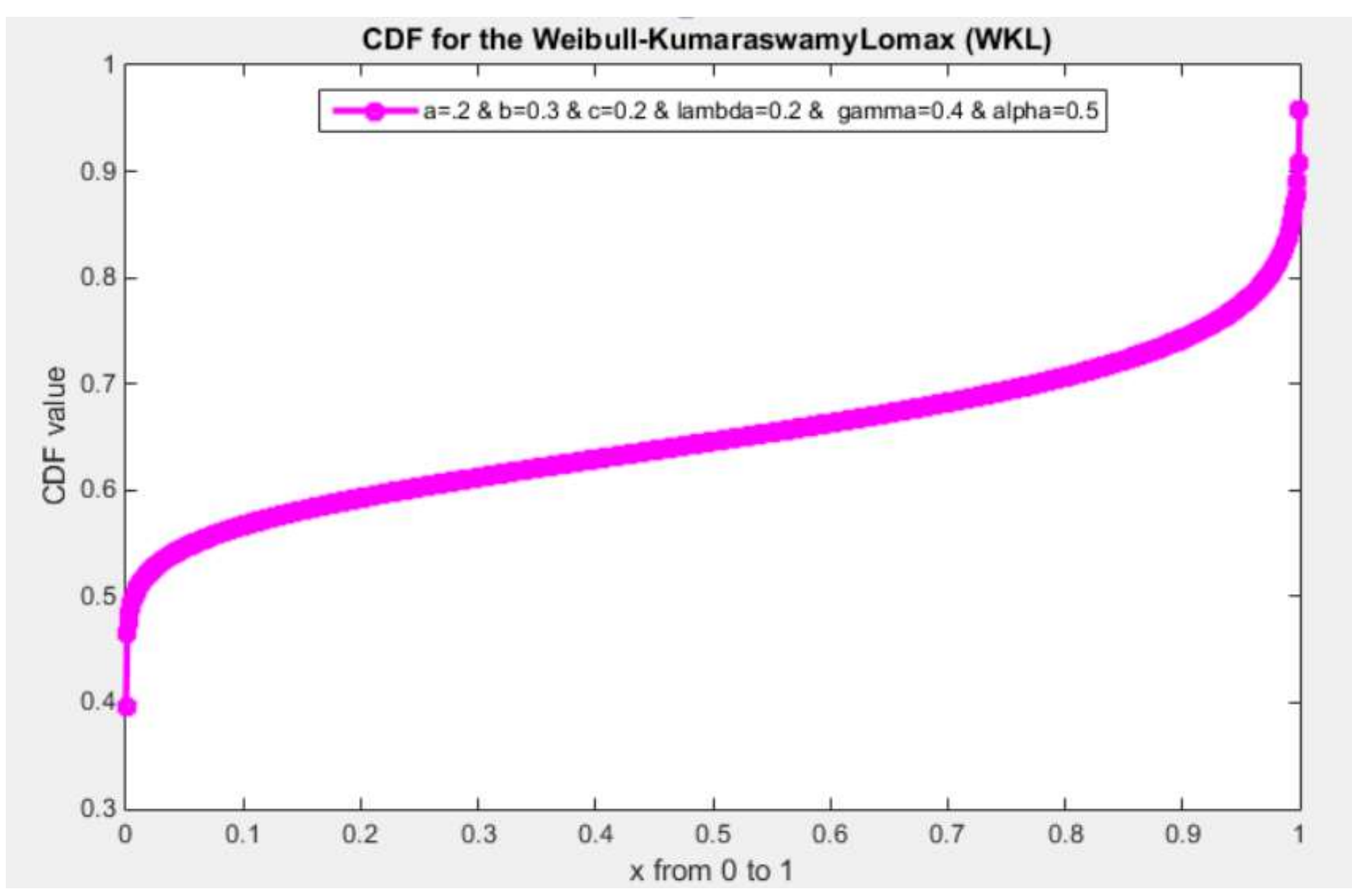

Fig. 25 shows the CDF when the parameters are a=0.2, b=0.3, c=0.2, lambda=0.2, gamma=0.4 and alpha=0.5

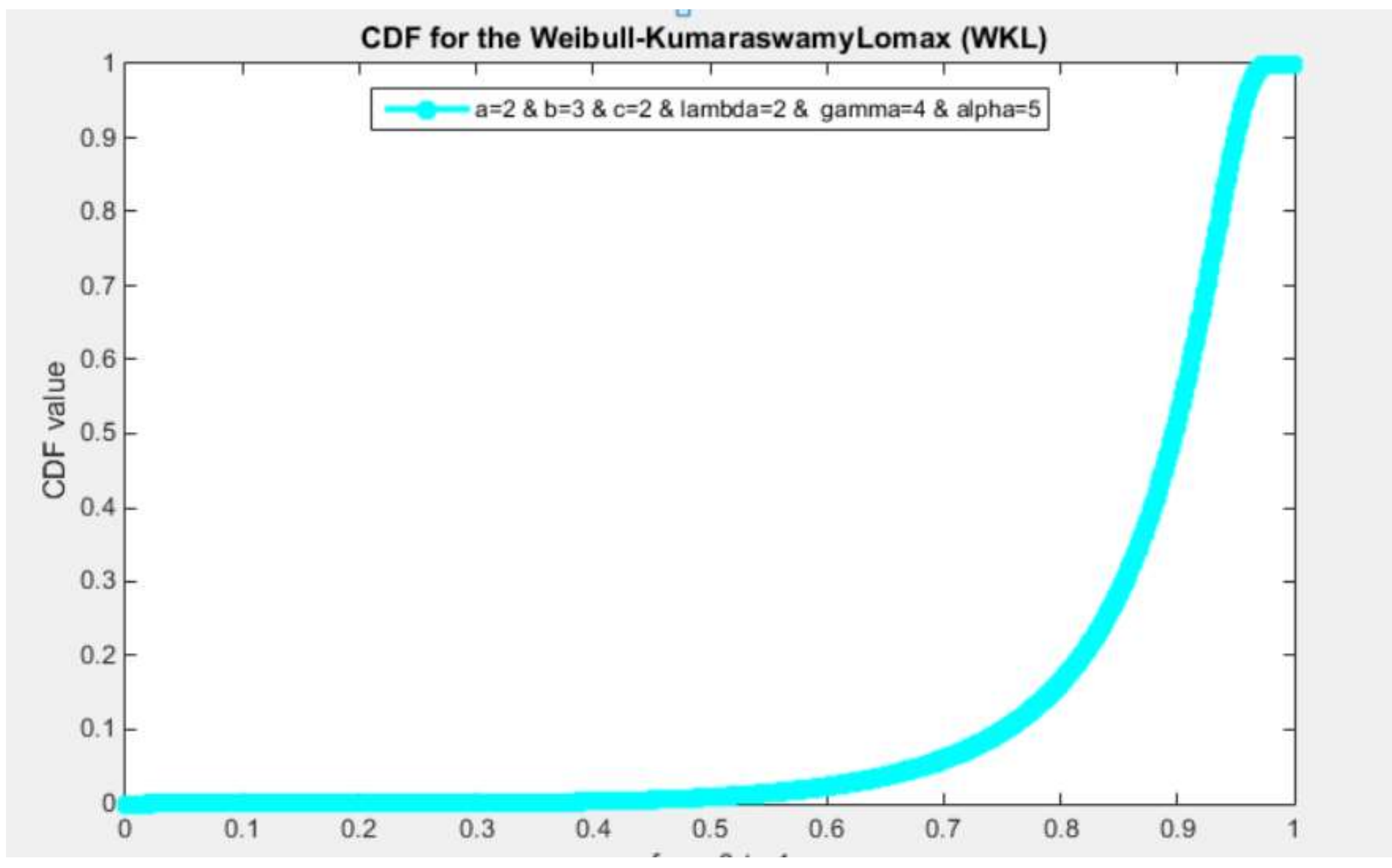

Fig. 26 shows the CDF when the parameters are a=2, b=3, c=2, lambda=2, gamma=4 and alpha=5

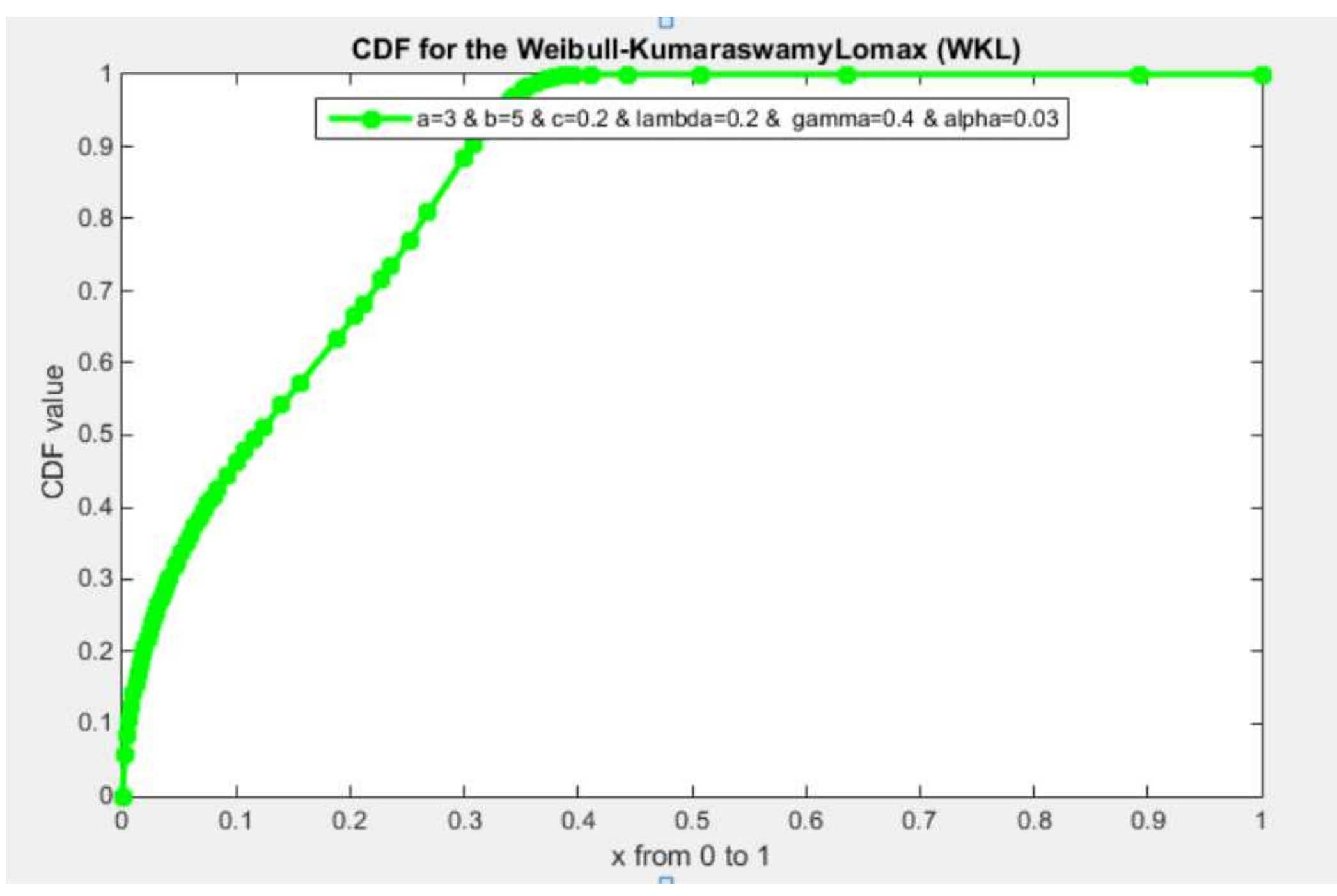

Fig. 27 shows the CDF when the parameters are a=3, b=5, c=0.2, lambda=0.2, gamma=0.4 and alpha=0.03

## 5. Results of Real Data Analysis and Discussion

The following is 56 observations of the death rate in Canada for the COVID 19 patients: 0.1622, 0.1159, 0.1897, 0.1260, 0.3025, 0.2190, 0.2075, 0.2241, 0.2163, 0.1262, 0.1627, 0.2591, 0.1989, 0.3053, 0.2170, 0.2241, 0.2174, 0.2541, 0.1997, 0.3333, 0.2594, 0.223, 0.2290, 0.1536, 0.2024, 0.2931, 0.2739, 0.2607, 0.2736, 0.2323, 0.1563, 0.2677, 0.2181, 0.3019, 0.2136, 0.2281, 0.2346, 0.1888, 0.2729, 0.2162, 0.2746, 0.2936, 0.3259, 0.2242, 0.1810, 0.2679, 0.2296, 0.2992, 0.2464, 0.2576, 0.2338, 0.1499, 0.2075, 0.1834, 0.3347, 0.2362. These observations have a mean value of 0.2305, a standard deviation value of 0.0520, a skewness value of -0.0873, a kurtosis value of 2.6537, a minimum value at 0.1159, and a maximum value of 0.3347. The 25$^{th}$ quantile is 0.2011, the median is 0.2262, and the 75$^{th}$ quantile is 0.2678. The author used fmincon function in MATLAB to estimate the parameters. Table (1) shows those results:

Table (1): Parameters estimators and their standard error and confidence interval

| parameters | estimate | Standard error | Lower 95 %CI | Upper 95% CI |
|---|---|---|---|---|
| $(a)$ | 6.8741 | 0.1234 | 6.6322 | 7.116 |
| $(b)$ | 23.876 | 4.4967 | 15.063 | 32.69 |
| $(c)$ | 0.73187 | 0.07677 | 0.58139 | 0.88236 |
| $(\lambda)$ | 3.2528 | 3.7132 | 0 | 10.531 |
| $(\alpha)$ | 3.028 | 0.41389 | 2.2167 | 3.8392 |
| $(\gamma)$ | 0.0019099 | 0.001109 | 0 | 0.00408 |

The final negative log-likelihood is -86.6007, the AIC is -161.2014, corrected AIC is -159.4871, the BIC is -149.0493, the HQIC is -159.4900. The KS test cannot reject that the WKL distribution fits the data with a p-value 0.5586. the variance covariance matrix of the estimated parameters is

$$var-cov=\begin{bmatrix} 0.0152 & 0.5499 & 0.0023 & -0.4515 & 0.0506 & 0.00126 \\ 0.5499 & 20.2203 & 0.087 & -16.6331 & 1.8595 & 0.0046 \\ 0.0023 & 0.0870 & 0.0059 & -0.0777 & 0.0077 & 0.000028 \\ -0.4515 & -16.6331 & -0.0777 & 13.7878 & -1.5240 & -0.0037 \\ 0.0506 & 1.8595 & 0.0077 & -1.5240 & 0.1713 & 0.00043 \\ 0.0012 & 0.0046 & 0.000028 & -0.0037 & 0.0004 & 0.0000012 \end{bmatrix}$$

The following figures (28-29) show the fitted PDF and the fitted CDF of the WKL distribution. The fitted curves fit the data well. This new distribution was compared with the Beta, Kumaraswamy, Topp –Loene distribution, Unit Lindley and Generalized Odd Median Based Unit Rayleigh type I (GOMBUR type I) as regards fitting the data. The estimators of the $a$ parameter for the Beta distribution was 14.5128 and for the $b$ parameter was 48.4899 with variance value for the $a$ parameter equals to 7.839 and for the $b$ parameter equals 100.6504. The values of the AIC, corrected –AIC, BIC and HQIC were respectively recording -167.88, -167.6536, -163.8293, -1663096. The negative LL was -85.94. The KS-test failed to reject that the distribution can fit the data well with p-value equals 0.6802. The covariance between the two shape parameters of the Beta distribution was estimated to be 27.531 which is high. After fitting the Kumaraswamy distribution to the data, the estimated $a$ parameter was 5.0309 with estimated variance equals 0.2719 while the estimated $b$ parameter was 1049.6 with estimated variance equals 523950. The estimated covariance between the two parameters was 370.2768 which is also high. The AIC, corrected-AIC, BIC, HQIC were respectively estimated to be -169.2, -168.9736, -165.1493, -167.6296, and the NLL was -86.6 and the Ks-test failed to reject the null hypothesis that the kumaraswamy can fit the data well with a p-value equals 0.5583. The Topp-Loene distribution and Unit Lindley distribution did not fit the data as KS-test for both of them reject the null hypothesis that each of them can fit the data with p-values equal to zero for both of the distributions. AIC, corrected-AIC, BIC, HQIC for the Topp-Leone were -46.3748, -46.3008, -443495, -45.5896 and NLL was 24.1874. As regards, Unit-Lindley the recorded AIC, corrected-AIC, BIC and HQIC were respectively -80.2707, -80.1966, -78.2453, -79.4855, and the NLL was 41.1353. The GOMBUR type I distribution, introduced by (Attia, 2025), fits the data well. The estimated n parameter equals 41.02961 with estimated variance equals 62.6379 and estimated $\alpha$ parameter of 1.4623 with estimated variance equals 0.0002376. The covariance between the two parameters was 0.0083 which is far less than the covariance between the two parameters of the Beta and the Kumaraswamy distributions. The NLL was -85.783 and the KS-test failed to reject the null hypothesis that the GOMBUR type I can fit the data well.

The new six-parameter distribution nearly gives comparable result to both Kumaraswamy and GOMBUR type I distribution although less than both of them. The AIC, corrected –AIC, BIC and HQIC for the new distribution is less than both of those two distributions but the NLL is the -68.6007. The variance-covariance matrix may show some correlation between some parameters but most of the covariance are less than that reported from fitting the beta and kumaraswamy distributions. And this is one of the advantages of adding new parameters to the baseline distribution and by generalizing this baseline distribution using different methodologies or techniques. Also adding new parameter may decrease the variance of some parameters and this improve the construction of the confidence interval. The determinant of the variance-covariance matrix of the fitted Beta distribution is 31.0432, the determinant of the fitted Kumaraswamy distribution is 5348.9 while the determinant for the fitted GOMBUR type I is 0.0148 and for the new distribution WKL distribution is 1.8519e-18 which is almost zero and this indicates there may be some correlation between the parameters.

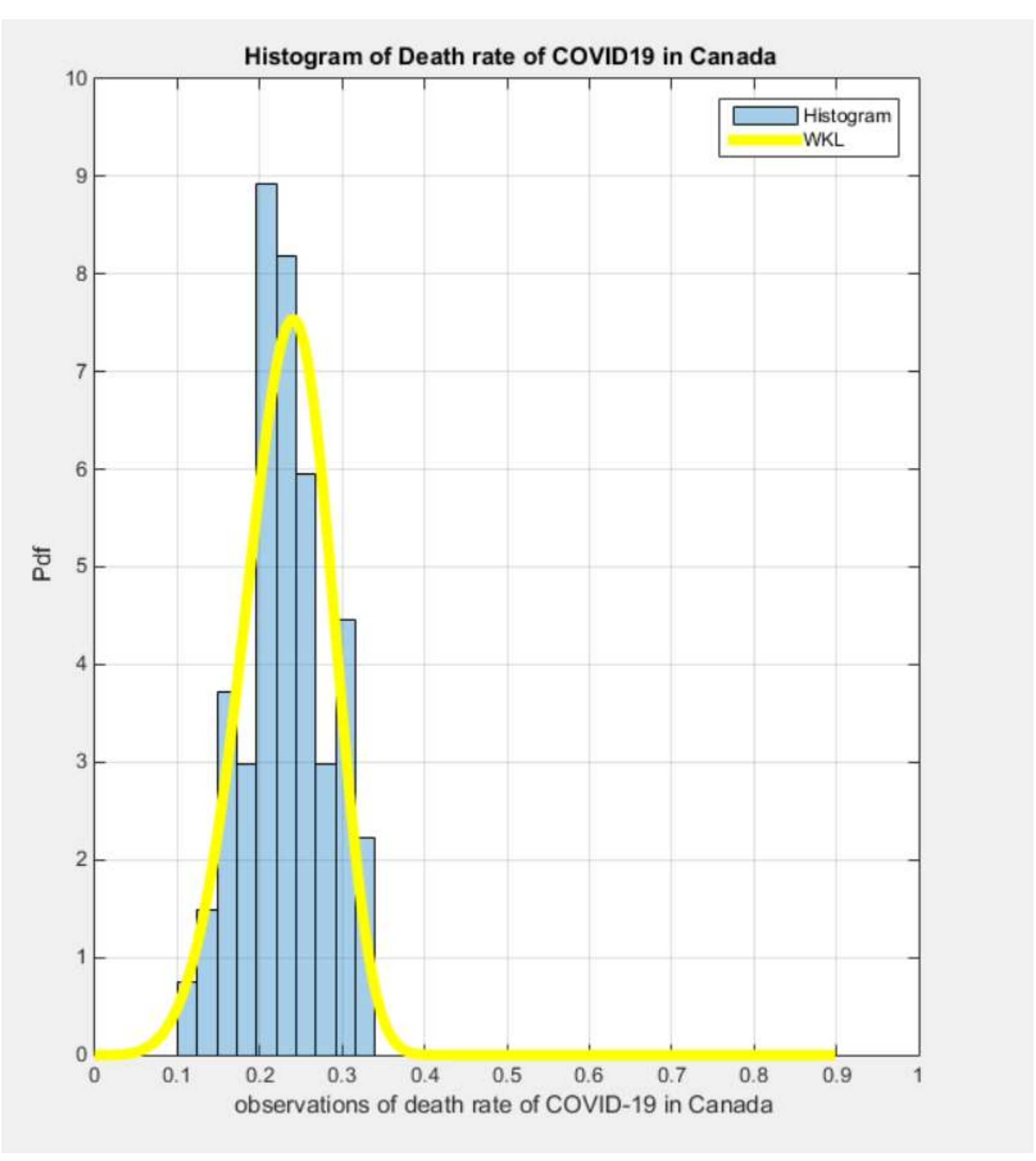

Fig. 28 shows that the data exhibit slight negative skewness. It also shows that the fitted PDF and the distribution fits the data well.

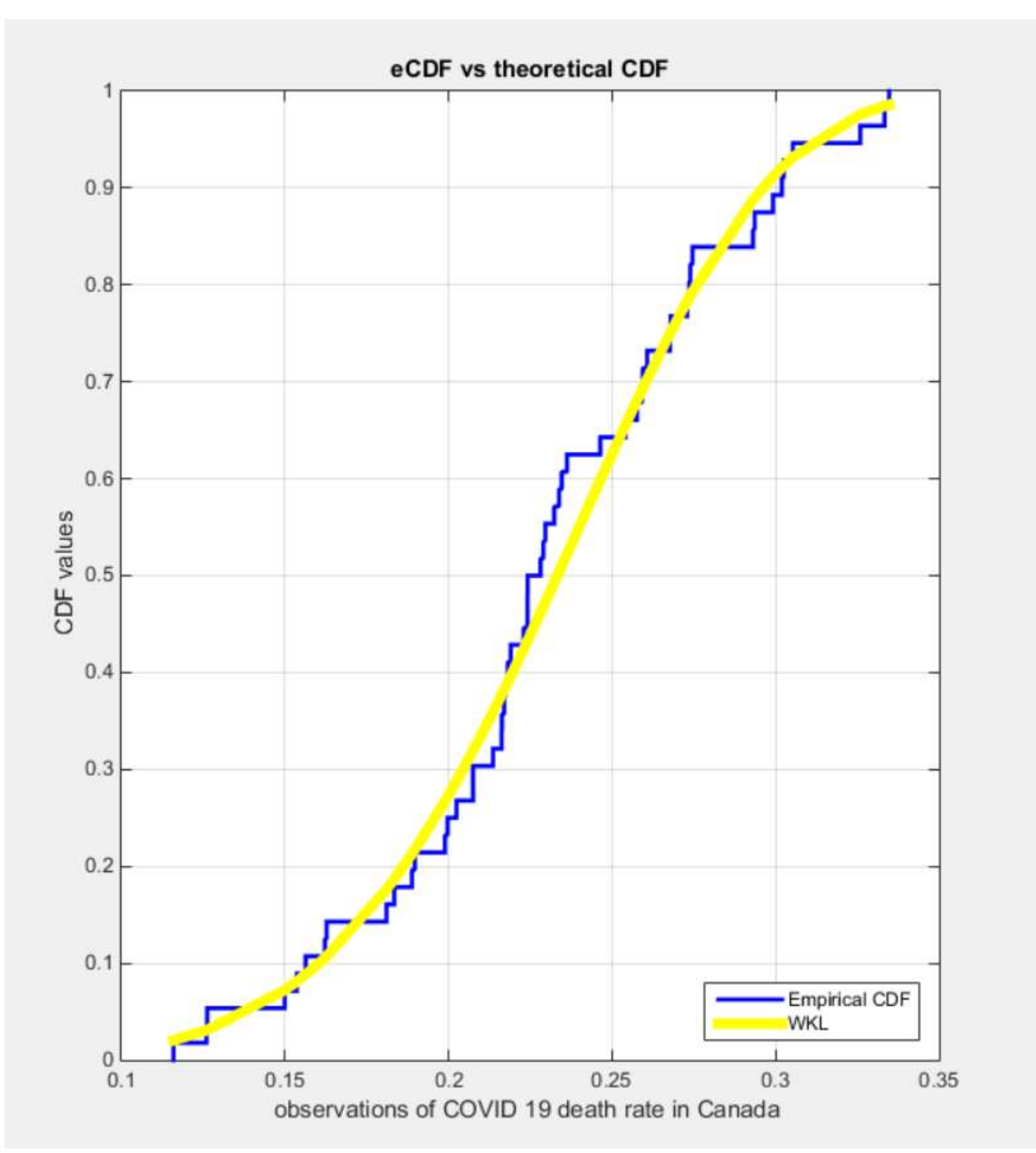

Fig. 29 shows that the fitted CDF and the distribution fits the data well.

## 6-Conclusion

The six parameter distribution has its parameter from the successive addition of the parameters of the constituent distributions, two shape parameters from the baseline distribution the kumaraswamy, the scale and shape parameters from the Lomax distribution, the link distribution of the quantile function, and the scale and shape parameters of the generator distribution, the Weibull distribution. These allow the curve of the PDF to fit and bend with the shape of the data and hence fit the data well as was shown from the Figures of the shape of the PDF. The addition of these new parameters decreases the covariance between many of the parameters although some parameters may exhibit high correlation. Also the addition of the parameters enhances the indices of validation like the AIC, corrected-AIC, BIC and HQIC. It also enhances the Log-Likelihood. In some cases, the baseline distribution may not fit the data but after generalization and addition of new parameters the new generalized distribution fit the data well. But in this case study and real data analysis, both the baseline and the generalized form fit the data well. And the baseline distribution has better indices than the generalized form. And this may raise a question is the new parameter addition may harbor some correlation between too many parameters and this need a solution for this correlated parameters. Addition of new parameters may decrease the determinant of the estimated variance –covariance matrix of the estimated parameters and may decrease the estimated variance of the estimated parameters thus enhancing the confidence interval construction. Analysis of real data reveals that the determinant of the estimated variance-covariance matrix of this new WKL distribution is nearly zero and this indicates that there is high possibility of correlation between some parameters and this requires further treatment of this singularity matrix and parameters correlations. Also, in this work the author introduces new approach for the linear series expansion for the PDF of the WKL to easy and simplify the statistical computations for calculating the raw moments and other properties of the distribution.

## 7-Future Work

Use other link functions or quantile distributions for the same baseline distribution the kumaraswamy and use different generator than the Weibull to compare those distributions and to see if the addition worth to be established. These comparisons enhance our understanding of the pros and cons of adding new parameters. Also simulation study to figure the behavior of the parameters.

**Declarations:**
**The author declares that the author has not used Artificial Intelligence tools in the creation of this manuscript.**
**Ethics approval and consent to participate**
Not applicable.
**Consent for publication**
Not applicable
**Data Availability Statement**
Data sharing not applicable to this article as no datasets were generated or analyzed during the current study.
**Competing interests**
The author declares no competing interests of any type.
**Funding**
No funding was received for this manuscript.
**Authors' contribution**
AI (Attia Iman) carried the conceptualization by formulating the goals, aims of the research article, formal analysis by applying the statistical, mathematical and computational techniques to synthesize and analyze the hypothetical data, carried the methodology by creating the model, software programming and implementation, supervision, writing, drafting, editing, preparation, and creation of the presenting work.
**Acknowledgement**
Not applicable